\documentclass{article}

\usepackage{amsmath}
\usepackage{latexsym}
\usepackage{graphicx}
\usepackage{psfrag}
\usepackage{amssymb}
\usepackage{amscd}

\def\<{\langle}
\def\>{\rangle}
\def\0{{{\bf 0}}}
\def\CD{{\mathcal D}}
\def\OO{{\mathcal O}}
\def\CE{{\mathcal E}}
\def\CF{{\mathcal F}}
\def\CL{{\mathcal L}}
\def\bCD{{\overline {\mathcal D}}}

\def\CA{{\mathcal A}}
\def\hCA{{\hat {\mathcal A}}}
\def\hCB{{\hat {\mathcal B}}}
\def\hCE{{\hat {\mathcal E}}}
\def\CB{{\mathcal B}}
\def\CM{{\mathcal M}}
\def\tCM{{\tilde {\mathcal M}}}

\def\CV{{\mathcal V}}

\def\AA{{\mathbb A}}
\def\CC{{\mathbb C}}

\def\FF{{\mathbb F}}

\def\QQ{{\mathbb Q}}
\def\RR{{\mathbb R}}
\def\TT{{\mathbb T}}
\def\ZZ{{\mathbb Z}}
\def\mm{{\mathfrak m}}
\def\pp{{\mathfrak p}}

\def\te{{\tilde e}}

\def\tX{{\tilde X}}
\def\tW{{\tilde W}}
\def\tU{{\tilde U}}
\def\hX{{\hat X}}

\def\tA{{\tilde A}}

\def\tlambda{{\tilde \lambda}}
\def\trho{{\tilde \rho}}

\newcommand{\spec}{\operatorname{Spec}}
\newcommand{\id}{\operatorname{id}}

\newcommand{\Fil}{\operatorname{Fil}}
\newcommand{\Hom}{\operatorname{Hom}}

\newcommand{\rank}{\operatorname{rank}}
\newcommand{\End}{\operatorname{End}}
\newcommand{\Fr}{\operatorname{Fr}}
\newcommand{\Fabs}{F_{\mbox{\rm \tiny abs}}}
\newcommand{\Lie}{\operatorname{Lie}}
\newcommand{\Frob}{\operatorname{Frob}}
\newcommand{\Ver}{\operatorname{Ver}}
\newcommand{\gal}{\operatorname{Gal}}
\newcommand{\Aut}{\operatorname{Aut}}
\newcommand{\DR}{\operatorname{DR}}
\newcommand{\GL}{\operatorname{GL}}
\newcommand{\GU}{\operatorname{GU}}

\newcommand{\red}{\mbox{\rm \tiny red}}

\newcommand{\cris}{\mbox{\rm \tiny cris}}
\newcommand{\et}{\mbox{\rm \tiny \'et}}

	{\begin{list}
		{\noindent\makebox[0mm][r]{\arabic{enumi}.}}
		{\usecounter{enumi} \topsep=1.5mm \itemsep=0mm}
	}
	{\end{list}}

	{\begin{list}
		{\noindent\makebox[0mm][r]{\arabic{enumi}.}}
		{\usecounter{enumi} \topsep=1.5mm \itemsep=-1mm}
	}
	{\end{list}}

\begin{document}

  \newtheorem{theorem}{Theorem}[section]
  \newtheorem{definition}[theorem]{Definition}
  \newtheorem{lemma}[theorem]{Lemma}
  \newtheorem{corollary}[theorem]{Corollary}
  \newtheorem{proposition}[theorem]{Proposition}
  \newtheorem{conjecture}[theorem]{Conjecture}
  \newtheorem{question}[theorem]{Question}
  \newtheorem{problem}[theorem]{Problem}

  \newtheorem{example}[theorem]{Example}
  \newtheorem{remark}[theorem]{Remark}

\newenvironment{proof}{{{\noindent{\it Proof.\ }}}}{{\hfill $\Box$}}

\title{Towards a geometric Jacquet-Langlands correspondence for unitary Shimura
varieties}  

\author{David Helm
\footnote{University of Texas at Austin; dhelm@math.utexas.edu}}

\maketitle

Let $G$ be a unitary group over a totally real field, and $X$ a Shimura
variety associated to $G$.  For certain primes $p$ of good reduction
for $X$, we construct cycles $X_{\tau_0,i}$ on the characteristic $p$
fiber of $X$.  These cycles are defined as the loci on which the
Verschiebung map has small rank on particular pieces of the Lie
algebra of the universal abelian variety on $X$.

The geometry of these cycles turns out to be closely related to Shimura
varieties for a {\em different} unitary group $G^{\prime}$, which is
isomorphic to $G$ at all finite places but not isomorphic to $G$ at
archimedean places.  More precisely, each cycle $X_{\tau_0,i}$ has
a natural desingularization $\tX_{\tau_0,i}$, which is ``almost''
isomorphic to a scheme parametrizing certain subbundles of the
Lie algebra of the universal abelian variety over a Shimura variety 
$X^{\prime}$ associated to $G^{\prime}$.

We exploit this relationship to construct an injection of the {\'e}tale
cohomology of $X^{\prime}$ into that of $X$.  This yields a 
geometric construction of ``Jacquet-Langlands transfers'' of
automorphic representations of $G^{\prime}$ to automorphic representations
of $G$.

\noindent
2000 MSC Classification: 11G18

\section{Introduction}

Suppose $G$ and $G^{\prime}$ are two algebraic groups over $\QQ$, isomorphic
at all finite places of $\QQ$ but not necessarily isomorphic at infinity.
The Jacquet-Langlands correspondence predicts, in many cases, that there
is a natural bijection between the automorphic representations of $G$ 
and those for $G^{\prime}$, such that if $\pi^{\prime}$ is the 
representation of $G^{\prime}$ corresponding to a representation $\pi$ of 
$G$, then $\pi_v$ is isomorphic to $\pi^{\prime}_v$ for all finite places $v$.

This correspondence is proven in many cases by comparing the trace formulas
for $G$ and $G^{\prime}$.  In this way one can conclude that there is
an isomorphism between suitable spaces of automorphic forms for $G$ and
$G^{\prime}$ as abstract representations, but not in any canonical fashion.
One might therefore hope for a more natural way of understanding this
correspondence.

For $\GL_2$, Serre \cite{Se} was the first to suggest an alternative approach
to Jacquet-Langlands, in the context of modular forms mod $p$.  By
considering modular forms as sections of a line bundle on a modular
curve, and restricting these sections to the supersingular locus
of this curve in characteristic $p$, Serre relates the space
of modular forms mod $p$ to a space of functions on this supersingular locus;
the latter can be interpreted as a space of ``algebraic modular forms''
for the quaternion algebra $B$ ramified at $p$ and infinity.  Ghitza has
since adapted this approach for symplectic groups \cite{Gh1},~\cite{Gh2}.

In contrast to the traditional approach to Jacquet-Langlands, the approach
of Serre and Ghitza yields canonical isomorphisms between spaces that
arise naturally from the geometry of Shimura varieties attached to the groups
under consideration (rather than simply a bijection of isomorphism classes
of representations.)

Another approach can be found in the  work of 
Ribet (\cite{level},~\cite{bimodules}).  Ribet finds a relationship between
the reductions at various primes of two Shimura curves associated to
two {\em different} quaternion algebras over $\QQ$.  He uses this to
obtain an explicit isomorphism between certain Hecke modules for the
two quaternion algebras, and thereby proves the Jacquet-Langlands 
correspondence in that setting.  This sharpening of the Jacquet-Langlands
correspondence is a key ingredient in his proof of level-lowering for
classical modular forms.

More recently, work of the author in~\cite{u2} adapts Ribet's techniques
to the case of a unitary group $G$ isomorphic 
(up to a factor of $\RR^{\times}$) to
a product of $U(1,1)$'s at infinity.  As with Ribet's approach, the key
is to understand the reduction of a Shimura variety $X$ attached to $G$
that has an analogue of $\Gamma_0(p)$ level structure at $p$.  We obtain
an explicit description of the global structure of the special fiber
in this setting: the irreducible components each are (nearly) isomorphic
to products of projective bundles over Shimura varieties $X^{\prime}$
attached to unitary groups $G^{\prime}$ that are isomorphic to $G$
at all finite places but not in general isomorphic to $G$ at infinity.
Via the theory of vanishing cycles, one can then relate the {\'e}tale
cohomology of $X$ to the {\'e}tale cohomology of the various $X^{\prime}$
that arise; the upshot is that the highest weight quotient of the {\'e}tale
cohomology of $X$ can be interpreted in terms of a space of algebraic
modular forms (over $\QQ_{\ell}$) for a unitary group $G^{\prime}$ that is
compact at infinity.

In this paper we present a different approach, that works for arbitrary
unitary groups, but proceeds by considering Shimura varieties at primes
of {\em good} reduction.  Given a Shimura variety $X$ arising from a unitary
group $G$, and a suitable prime $p$ of good reduction, we consider cycles
$X_{\tau_0,i}$ on the characteristic $p$ fiber of $X$.  (Here $i$ is a
positive integer and $\tau_0$ determines a map 
$p_{\tau_0}: \OO_F \rightarrow \overline{\FF}_p$,
where $F$ is the $CM$-field arising in the definition of $G$.)  Loosely
speaking, $X_{\tau_0,i}$ is the locus of abelian varieties $A$ (with 
$\OO_F$-action) such that the space $\Hom(\alpha_p,A[p])_{p_{\tau_0}}$
of maps on which $\OO_F$ acts via $p_{\tau_0}$ has dimension at least 
$i$ larger than the ``expected dimension''.  Alternatively,
$X_{\tau_0,i}$ can be thought of as the locus of abelian varieties $A$
such that 
$$\Ver: \Lie(A^{(p)})_{p_{\tau_0}} \rightarrow \Lie(A)_{p_{\tau_0}}$$
has rank at least $i$ less than the ``expected rank''.  Such cycles
are closed strata in the so-called ``a-number'' stratification of $X$,
and their local structure has been studied extensively.

Our approach requires an understanding of the {\em global} structure
of these strata in addition to the local structure.  Questions of
this nature are not nearly as well-understood; fortunately in the
cases we are interested in they can be attacked by fairly standard
techniques.  We construct a natural desingularization of each
cycles $\tX_{\tau_0,i}$.
The geometry of this desingularization turns out to be
closely related to the geometry of a Shimura variety $X^{\prime}$
arising from a different unitary group $G^{\prime}$.  
As in~\cite{u2}, $G^{\prime}$ is isomorphic to $G$ at finite
places but not at infinity.  In particular, we construct a scheme
$(X^{\prime})^{\tau_0,i}$, defined naturally in terms of the universal
abelian variety over $X^{\prime}$, such that there exists a scheme
$\hX_{\tau_0,i}$, together with maps: 
$$\hX_{\tau_0,i} \rightarrow \tX_{\tau_0,i}$$
$$\hX_{\tau_0,i} \rightarrow (X^{\prime})^{\tau_0,i}$$
that are bijections on points and isomorphisms on {\'e}tale cohomology.
Loosely speaking, this says that $\tX_{\tau_0,i}$ and $(X^{\prime})^{\tau_0,i}$
are ``nearly'' isomorphic.  The
fibers of $(X^{\prime})^{\tau_0,i}$ over $X^{\prime}$ are Grassmannians
of various dimensions.  (This generalizes results of~\cite{u2} for the
case of $U(2)$ Shimura varieties).
 
Rather attempting to establish a geometric Jacquet-Langlands
correspondence as in~\cite{u2}, by way of a suitable ``Deligne-Rapoport 
model'', we use the above construction to give an explicit injection
of the {\'e}tale cohomology of $X^{\prime}$ into that of $X$. 
(Theorem~\ref{thm:main}.) 
The existence of this map follows from a general construction in {\'e}tale
cohomology; its injectivity is more difficult to prove.  The key
ingredients are the Leray spectral sequence for $(X^{\prime})^{\tau_0,i}$
and the Thom-Porteus formula, which allows us to compute the
self-intersection of $X_{\tau_0,i}$.  This argument is the main focus
of sections~\ref{sec:cohomology} and~\ref{sec:porteus}.

We obtain cases of Jacquet-Langlands transfer as an easy corollary of
the existence of this injection. (Theorem~\ref{thm:JL})  In particular
we show that for every cohomological automorphic representation 
$\pi^{\prime}$ of $G^{\prime}$, there is an automorphic representation
$\pi$ of $G$ such that $\pi$ and $\pi^{\prime}$ agree at all finite places.
Our approach suffers from the limitation that it is only possible to
transfer such representations from one Shimura variety to another
Shimura variety of higher dimension; going in the other direction
requires some way of controlling the {\em image} of the map constructed in
Theorem~\ref{thm:main}, which we do not yet have at our disposal. 

This approach appears to cover some new cases of Jacquet-Langlands transfer
that have not yet appeared in the literature.  In particular the traditional
trace formula approach to Jacquet-Langlands runs into difficulty with
unitary groups that are not compact at infinity. On the other hand,
Harris and Labesse (\cite{HL}, particularly Theorems 2.1.2, 3.1.6, and 
Proposition 3.1.7) have used base change techniques 
to establish rather general Jacquet-Langlands results for unitary groups,
but need the representation under consideration to be supercuspidal at 
certain places.

It should also be noted that whereas traditional approaches to
Jacquet-Langlands yield a bijection of isomorphism classes of
automorphic representations, our approach yields information
about a natural map between spaces that arise naturally in geometry,
and have considerable arithmetic interest.  One might therefore hope
for arithmetic applications of this approach, in analogy with
the application of Ribet's results on character groups to
level-lowering.

\section{Basic definitions and properties} \label{sec:basic}

We begin with the definition and basic properties of unitary
Shimura varieties.

Fix a totally real field $F^+$, of degree $d$ over $\QQ$.  Let $E$ be an
imaginary quadratic extension of $\QQ$, of discriminant $D$, and let $x$ be
a square root of $D$ in $E$.  Let $F$ be the field $EF^+$.

Also fix a square root $\sqrt{D}$ of $D$ in $\CC$.  Then any embedding
$\tau: F^+ \rightarrow \RR$ induces two embeddings 
$p_{\tau}, q_{\tau}: F \rightarrow \CC$, via
\begin{eqnarray*}
p_{\tau}(a+bx) & = & \tau(a) + \tau(b)\sqrt{D} \\
q_{\tau}(a+bx) & = & \tau(a) - \tau(b)\sqrt{D}.
\end{eqnarray*}

Fix an integer $n$, and an $n$-dimensional $F$-vector space $\CV$, 
equipped with an alternating, nondegenerate
pairing $\langle, \rangle: \CV \times \CV \rightarrow \QQ.$  We require that
$$\langle \alpha x,y\rangle = \langle x, \overline{\alpha} y \rangle$$ 
for all $\alpha$ in $F$.

Each embedding $\tau: F^+ \rightarrow \RR$ gives us a complex vector
space $\CV_{\tau} = \CV \otimes_{F^+,\tau} \RR$.  The pairing $\langle,\rangle$
on $\CV$ is the ``imaginary part'' of a unique Hermitian pairing 
$[,]_{\tau}$ on $\CV_{\tau}$; we denote the
number of $1$'s in the signature of $[,]_{\tau}$ by $r_{\tau}(\CV)$, and the
number of $-1$'s by $s_{\tau}(\CV)$.  If $\CV$ is obvious from the context,
we will often omit it, and denote $r_{\tau}(\CV)$ and $s_{\tau}(\CV)$
by $r_{\tau}$ and $s_{\tau}$.  We fix a $\hat \OO_F$-lattice $T$
inside $\CV(\AA^f_{\QQ})$, such that $\lambda$ induces a map
$T \rightarrow \Hom(T, \hat \ZZ)$.

Let $G$ be the algebraic group over $\QQ$ such that for any $\QQ$-algebra
$R$, $G(R)$ is the subgroup of $\Aut_F(\CV \otimes_{\QQ} R)$ consisting
of all $g$ such that there exists an $r$ in $R^{\times}$ with
$\langle gx,gy \rangle = r\langle x,y \rangle$ for all $x$ and $y$ in
$V \otimes_{\QQ} R$.  The discussion in the previous paragraph shows
that $G(\RR)$ is the subgroup of 
$$\prod_{\tau: F^+ \rightarrow \RR} \GU(r_{\tau}, s_{\tau})$$
consisting of those tuples $(g_{\tau})_{\tau: F^+ \rightarrow \RR}$
such that the ``similitude ratio'' of $g_{\tau}$ is the same for all
$\tau$.

Now fix a compact open subgroup $U$ of $G(\AA^f_{\QQ})$, preserving $T$,
and consider the Shimura variety associated to $G$ and $U$.  If $U$ is
sufficiently small, this variety can be thought of as a fine moduli
space for abelian varieties with PEL structures.  We now describe such
a model over a suitable ring of Witt vectors. 

Fix a prime $p$ split in $E$, such that the cokernel of the map
$T \rightarrow \Hom(T, \hat \ZZ)$ is supported away from $p$ and
such that $U_p$ is equal to the subgroup of all elements of $G(\QQ_p)$
preserving $T(\QQ_p)$.  Also fix a finite field $k_0$ of characteristic
$p$ large enough to contain subfields isomorphic to each of the
residue fields of $\OO_F/p$, and an identification of the Witt
vectors $W(k_0)$ with a subring of $\CC$.  This choice of identification
induces a bijection of the set of archimedean places of $F$ with
the set of algebra maps $\OO_F \rightarrow W(k_0)$.  In an abuse
of notation we will use the symbols $p_{\tau}$ and $q_{\tau}$ to
represent both the embeddings $F \rightarrow \CC$ defined above, and
the maps $\OO_F \rightarrow W(k_0)$ defined above.  

Note that if $S$ is a $W(k_0)$-scheme, and $M$ is a $W(k_0)$-module
with an action of $\OO_F$, then $M$ splits as a direct sum
$$M = \bigoplus_{\tau} M_{p_{\tau}} \oplus M_{q_{\tau}},$$
where $\OO_F$ acts of $M_{p_{\tau}}$ (resp. $M_{q_{\tau}}$)
by $p_{\tau}: \OO_F \rightarrow W(k_0)$ (resp. 
$q_{\tau}: \OO_F \rightarrow W(k_0)$.)

Consider the functor that associates to each $W(k_0)$-scheme $S$ the
set of isomorphism classes of triples $(A,\lambda, \rho)$ where:
\begin{enumerate}
\item $A$ is an abelian scheme over $S$ of dimension $nd$, with an action
of $\OO_F$
\item $\lambda$ is a polarization of $A$, of degree prime to $p$, such
that the Rosati involution associated to $\lambda$ induces complex
conjugation on $\OO_F \subset \End(A)$.
\item $\rho$ is a $U$-orbit of isomorphisms 
$T^{(p)} \rightarrow T_{\hat \ZZ^{(p)}} A$, sending the Weil pairing on 
$T_{\hat \ZZ^{(p)}} A$
to a scalar multiple of the pairing $\langle, \rangle$ on $T^{(p)}$.  (Here
$T_{\hat \ZZ^{(p)}} A$ denotes the product over all $l \neq p$ of the $l$-adic 
Tate modules of $A$, and the superscript $(p)$ denotes the non-pro-$p$
part of $T$ or $\hat \ZZ$.) 
\item When considered as an endomorphism of $\Lie(A/S)$, an element
$\alpha$ of $\OO_F$ has characteristic polynomial
$$\prod_{\tau: F^+ \rightarrow \RR} (x - p_{\tau}(\alpha))^{r_{\tau}(\CV)}
(x - q_{\tau}(\alpha))^{s^{\tau}(\CV)}.$$
\end{enumerate}

Since $\Lie(A/S)$ is an $\OO_F \otimes S$-module, we can reprhrase
the characteristic polynomial condition as follows:
For each $\tau: F^+ \rightarrow \RR$, we have
\begin{eqnarray*}
\rank_S \Lie(A/S)_{p_{\tau}} & = & r_{\tau}(\CV)\\
\rank_S \Lie(A/S)_{q_{\tau}} & = & s_{\tau}(\CV).
\end{eqnarray*}

If $U$ is sufficiently small, this functor is represented by a smooth
$W(k_0)$-scheme, which we denote $X_U(\CV)$.  It is a model
for the Shimura variety associated to $G$ and $U$, over $W(k_0)$.  
Henceforth we will refer to such an object as a 'unitary Shimura variety'.
Its dimension is given by the formula 
$$\dim X_U(\CV) = \sum_{\tau: F^+ \rightarrow \RR} r_{\tau}s_{\tau}.$$  

\begin{remark} \rm The scheme $X_U(\CV)$
depends not only on $U$ and $\CV$ but on all of the choices we have
made in this section.  To avoid clutter, we have chosen to suppress 
most of these choices in our notation.
\end{remark}
 
\section{Dieudonn{\'e} theory and points on $X_U(\CV)$} 
\label{sec:points}

Let $k$ be a perfect field containing $k_0$, and let $(A,\lambda,\rho)$
be a $k$-valued point of $X_U(\CV)$.  We begin by studying the
(contravariant) Dieudonn{\'e} modules of $A[p]$ and $A[p^{\infty}]$.

Let $\CD_A$ denote the contravariant Dieudonn{\'e} module of $A[p^{\infty}]$.
It is a free $W(k)$-module of rank $2nd$, equipped with endomorphisms
$F$ and $V$, that satisfy $FV = VF = p$.  These endomorphisms
do not commute with the action of $W(k)$, but instead satisfy:
\begin{eqnarray*}
F\alpha & = & \alpha^{\sigma}F\\
V\alpha^{\sigma} & = & \alpha V,
\end{eqnarray*}
where $\alpha \in W(k)$, and the superscript $\sigma$ denotes the Witt 
vector Frobenius.

The $\OO_F$-action on $A$ induces an $\OO_F$ action on $\CD_A$;
we therefore have a direct sum decomposition:
$$\CD_A = 
\bigoplus_{\tau: F^+ \rightarrow \RR} (\CD_A)_{p_{\tau}} \oplus 
(\CD_A)_{q_{\tau}}.$$

For $\tau: F^+ \rightarrow \RR$, let $p_{\sigma\tau}$ denote the map
$\OO_F \rightarrow W(k_0)$ obtained by taking
the map $\OO_F \rightarrow W(k_0)$ corresponding to $p_{\tau}$ and
composing it with the Witt vector Frobenius.  Define $q_{\sigma\tau}$
similarly.  Then the $\sigma$-linearity properties of $F$ and $V$ mean
that they induce maps:
\begin{eqnarray*}
F: (\CD_A)_{p_{\tau}} & \rightarrow & (\CD_A)_{p_{\sigma\tau}} \\
V: (\CD_A)_{p_{\sigma\tau}} & \rightarrow & (\CD_A)_{p_{\tau}},
\end{eqnarray*}
and similarly for the $q_{\tau}$.  Since $FV = VF = p$, we find that 
$(\CD_A)_{p_{\tau}}$ and $(\CD_A)_{p_{\sigma\tau}}$ have the same rank
for all $\tau$, as do $(\CD_A)_{q_{\tau}}$ and $(\CD_A)_{q_{\sigma\tau}}$.

If we fix a prime $\pp$ of $\OO_F$ over $p$, then the 
Dieudonn{\'e} module of $A[\pp^{\infty}]$ is the direct sum of 
$(\CD_A)_{p_{\tau}}$ for those $p_{\tau}$ (or possibly $q_{\tau}$) 
for which the preimage of the ideal $(p)$ of $W(k_0)$ under the
corresponding map $\OO_F \rightarrow W(k_0)$ is $\pp$.  These form
a single orbit under the action of $\sigma$ described above, so
they all have the same rank.  But since the height of
$A[\pp^{\infty}]$ is $n$ times the residue class degree of $\pp$ over
$p$, it follows that $(\CD_A)_{p_{\tau}}$ and $(\CD_A)_{q_{\tau}}$
are free $W(k)$-modules of rank $n$ for all $\tau$. 

Now consider the quotient $\bCD_A = \CD_A/p\CD_A$.  It is canonically
isomorphic to the Dieudonn{\'e} module of $A[p]$.  The above discussion
shows that for each $\tau$, $(\bCD_A)_{p_{\tau}}$ and $(\bCD_A)_{q_{\tau}}$
are $n$-dimensional $k$-vector spaces.  Moreover, Oda~\cite{oda} has shown
that there is a natural isomorphism $H^1_{\DR}(A/k) \cong \bCD_A$, and
that this isomorphism identifies the Hodge flag $\Lie(A/k)^* \subset
H^1_{\DR}(A/k)$ with the subspace $V\bCD_A$ of $\bCD_A$.

In particular, we have $\dim V((\bCD_A)_{p_{\sigma\tau}}) = 
\dim \Lie(A/k)^*_{p_{\tau}} = r_{\tau}$.  Since the image of $V$ is equal
to the kernel of $F$ on $\bCD_A$, we also have
$\dim F((\bCD_A)_{p_{\tau}}) = s_{\tau}.$

Thus, for each $\tau$, $(\bCD_A)_{p_{\tau}}$ is an $n$-dimensional
$k$-vector space with two distinguished subspaces, 
$F_{\tau} = F((\bCD_A)_{p_{\sigma^{-1}\tau}})$ (of dimension 
$s_{\sigma^{-1}\tau}$), and
$V_{\tau} = V((\bCD_A)_{p_{\sigma\tau}})$ (of dimension $r_{\tau}$.)

Fix a particular $\tau_0$, and assume, for the rest of the paper,
that $r_{\sigma^{-1}\tau_0} \leq r_{\tau_0}$.  (If this does not hold,
then $s_{\sigma^{-1}\tau_0} \leq s_{\tau_0}$, and everything that
follows will still be true once one reverses the roles of $p_{\tau}$
and $q_{\tau}$.)  In this case, if $F_{\tau_0}$ and $V_{\tau_0}$
are in general position with respect to each other, then their sum
will span all of $(\bCD_A)_{\tau_0}$.  Of course, $F_{\tau_0}$ and
$V_{\tau_0}$ need not be in general position with respect to one
another, which motivates the following definition:

\begin{definition} Let $i$ be an integer between $0$ and 
$\min(r_{\sigma^{-1}\tau},
s_{\tau})$, inclusive.  A point $(A,\lambda,\rho)$ is 
$(\tau_0,i)$-special if $\dim F_{\tau_0} + V_{\tau_0} \leq n-i$.  A
subspace $H$ of $(\bCD_A)_{p_{\tau_0}}$ is $(\tau_0,i)$-special
if it has dimension $n-i$ and contains both $F_{\tau_0}$ and $V_{\tau_0}$.
\end{definition}

Note that $(A,\lambda,\rho)$ admits an $H$ that is $(\tau_0,i)$-special
if and only if $(A,\lambda,\rho)$ itself is $(\tau_0,i)$-special, and
that such an $H$ will be unique if and only if $(A,\lambda,\rho)$ is
$(\tau_0,i)$-special but not $(\tau_0,i+1)$-special.

Suppose we have $(A,\lambda,\rho)$, along with a $(\tau_0,i)$-special
$H$ for this abelian variety.  Define a submodule $M_H$ of $\bCD_A$
as follows:
\begin{enumerate}
\item $(M_H)_{p_{\tau_0}} = H$
\item $(M_H)_{p_{\tau}} = (\bCD_A)_{p_{\tau}}$ for $\tau \neq \tau_0$
\item $(M_H)_{q_{\tau}} = (M_H)_{p_{\tau}}^{\perp}$, where $\perp$ denotes
orthogonal complement under the perfect pairing
$(\bCD_A)_{p_{\tau}} \times (\bCD_A)_{q_{\tau}} \rightarrow k$ induced
by the polarization $\lambda$.
\end{enumerate}

It is clear that $M_H$ is stable under $W(k)$, $\OO_F$, $F$, and $V$.
In particular, it is a Dieudonn{\'e} submodule of $\bCD_A$.  We thus
obtain an exact sequence:
$$0 \rightarrow M_H \rightarrow \bCD_A \rightarrow \bCD_K \rightarrow 0$$
where $\bCD_K$ is the Dieudonn{\'e} module of a group scheme $K$ over $k$.
The surjection $\bCD_A \rightarrow \bCD_K$ corresponds to an inclusion
of $K$ in $A[p]$; henceforth we identify $K$ with its image in $A[p]$.
Since $M_H$ is a maximal isotropic subspace of $\bCD_A$ under the pairing
induces by $\lambda$, $K$ is a maximal isotropic subgroup of $A[p]$
(under the Weil pairing induced by $\lambda$).

Let $B = A/K$, and let $f: A \rightarrow B$ denote the quotient map.  
Since $K \subset A[p]$, multiplication by $p$ (considered as an
endomorphism of $A$) factors through $f$.  In this way we obtain
a map $f^{\prime}$ such that $ff^{\prime} = f^{\prime}f = p$.
Note that $f^{\prime}(B[p])$ is equal to $K$.

Consider the polarization $(f^{\prime})^{\vee}\lambda f^{\prime}$ of
$B$.  For any $\alpha, \beta$ in $B[p]$, we have
$$\langle \alpha, (f^{\prime})^{\vee}\lambda f^{\prime} \beta \rangle_B =
  \langle f^{\prime} \alpha, \lambda f^{\prime} \beta \rangle_A.$$
The right-hand side vanishes identically since $K$ is isotropic and
$f^{\prime}(B[p]) = K$.  Thus $B[p]$ lies in the kernel of
$(f^{\prime})^{\vee}\lambda f^{\prime}$, and so there is a unique
polarization $\lambda^{\prime}$ of $B$ such that
$p\lambda^{\prime} = (f^{\prime})^{\vee} \lambda f^{\prime}.$
(Note that $\lambda^{\prime}$ can also be characterized as the unique
polarization of $B$ such that $p\lambda = f^{\vee} \lambda^{\prime} f$.)
The degree of $\lambda$ is easily seen to be prime to $p$.

\begin{proposition} Suppose that $\sigma \tau_0 \neq \tau_0$.  Then
\begin{enumerate}
\item $\dim \Lie(B/k)_{p_{\tau_0}} = r_{\tau_0} + i$.
\item $\dim \Lie(B/k)_{p_{\sigma^{-1}\tau_0}} = r_{\sigma^{-1}\tau_0} - i$.
\item $\dim \Lie(B/k)_{p_{\tau}} = r_{\tau}$ for $\tau$ not equal to
$\tau_0$ or $\sigma^{-1}\tau_0$.
\item $\dim \Lie(B/k)_{q_{\tau}} = n - \dim \Lie(B/k)_{p_{\tau}}$ for 
all $\tau$.
\end{enumerate}
\end{proposition}

\begin{proof}
The quotient map $f: A \rightarrow B$ induces a map
$\bCD_B \rightarrow \bCD_A$, where $\bCD_B$ is the Dieudonn{\'e} module
of $B[p]$.  The image of this map is precisely $M_H$.  On the level
of $p$-divisible groups, therefore, $f$ induces an inclusion
of $\CD_B$ into $\CD_A$, that identifies $\CD_B$ with the submodule
of $\CD_A$ consisting of those elements whose images in $\bCD_A$ lie in
$M_V$.  We identify $\CD_B$ with this submodule for the remainder of
the argument.

By the isomorphism between Dieudonn{\'e} modules and DeRham cohomology, 
$$\dim \Lie(B/k)_{p_{\tau}} = 
  \dim V((\CD_B)_{p_{\sigma\tau}})/p(\CD_B)_{p_{\tau}}.$$
On the other hand, we have:
\begin{enumerate}
\item $(\CD_B)_{p_{\tau}} = (\CD_A)_{p_{\tau}}$ for $\tau \neq \tau_0$.
\item $(\CD_A)_{p_{\tau_0}}/(\CD_B)_{p_{\tau_0}}$ has dimension $i$.
\item $V((\CD_A)_{p_{\sigma\tau}})/p(\CD_A)_{p_{\tau}}$ has dimension
$r_{\tau}$ for all $\tau$.
\end{enumerate}

Statements (1), (2), and (3) of the proposition follow immediately
from the above paragraph.  Statement (4) follows from the existence of
the prime-to-$p$ polarization $\lambda^{\prime}$ on $B$.
\end{proof}

Note that if $\sigma \tau_0 = \tau_0$, then the result above fails.  (In
particular, the proof of the result shows in this case that 
$\dim \Lie(B/K)_{p_{\tau}} = r_{\tau}$ for all $\tau$.)  Since the
above proposition is crucial to our argument, we
assume, for the remainder of the paper, that $\sigma \tau_0 \neq \tau_0$.

The upshot of the above proposition is that $(B,\lambda^{\prime})$ is
``nearly'' a $k$-valued point a unitary Shimura variety.  It lacks
only a level structure.  We cannot define such a level structure in
terms of $\CV$, however, as $r_{\tau_0}(\CV) = r_{\tau_0}$ but
$\dim \Lie(B/k)_{p_{\tau_0}} = r_{\tau_0} + i$.  We thus invoke the following
lemma, proven in the appendix of~\cite{u2}:

\begin{lemma} \label{lemma:lierank}
There exists an $n$-dimensional $F$-vector space
$\CV^{\prime}$, together with a pairing $\langle, \rangle^{\prime}$
satisfying the conditions of section~\ref{sec:basic}, such that:
\begin{enumerate}
\item $r_{\tau_0}(\CV^{\prime}) = r_{\tau_0} + i$.
\item $r_{\sigma^{-1}\tau_0}(\CV^{\prime}) = r_{\sigma^{-1}\tau_0} - i$.
\item $r_{\tau}(\CV^{\prime}) = r_{\tau}$ for $\tau$ not equal to $\tau_0$
or $\sigma^{-1}\tau_0$.
\item There exists an isomorphism $\phi$ of $\CV(\AA^f_{\QQ})$ with
$\CV^{\prime}(\AA^f_{\QQ})$ that takes the pairing $\langle,\rangle$
to a scalar multiple of $\langle,\rangle^{\prime}$.
\end{enumerate}
\end{lemma}

We fix, once and for all, a $\CV^{\prime}$, $\langle,\rangle^{\prime}$
and $\phi$ as in the lemma.  Let $T^{\prime} = \phi(T)$, and let $G^{\prime}$
be the algebraic group such that for each $\QQ$-algebra $R$, 
$G^{\prime}(R)$ is the subset of $\Aut_F(\CV^{\prime} \otimes_{\QQ} R)$
consisting of those automorphisms that send $\langle,\rangle^{\prime}$
to a scalar multiple of itself.  Then $\phi$ induces an isomorphism 
$G(\AA^f_{\QQ}) \cong G^{\prime}(\AA^f_{\QQ})$, and this identifies
$U$ with a subgroup $U^{\prime}$ of $G^{\prime}$.  If $\rho$ is
a $U$-level structure on $(A,\lambda)$, then it follows from this construction
that $f \circ \rho \circ \phi^{-1}$ is a $U$-level structure on
$(B,\lambda^{\prime})$.  In particular, 
$(B,\lambda^{\prime}, f \circ \rho \circ \phi^{-1})$ is a $k$-valued point 
of the unitary Shimura variety $X_{U^{\prime}}$ associated to the subgroup
$U^{\prime}$ of $G^{\prime}$.

The map that associates to each $(A,\lambda,\rho,V)$ the point
$(B,\lambda^{\prime}, f \circ \rho \circ \phi^{-1})$ is not in general
a bijection.  We will now proceed to remedy this, by describing the
extra information needed to recover $(A,\lambda,\rho,V)$ from
$(B,\lambda^{\prime}, f \circ \rho \circ \phi^{-1})$. 

\begin{definition} Let $(B,\lambda^{\prime}, \rho^{\prime})$ be
a point on $X_{U^{\prime}}(k)$.  A subspace $W$ of $(\bCD_B)_{p_{\tau_0}}$
is {\em $(\tau_0,i)$-constrained} if it has dimension $i$ and
is contained in both $V((\bCD_B)_{p_{\sigma\tau_0}})$ and 
$F((\bCD_B)_{p_{\sigma^{-1}\tau_0}})$.
\end{definition}

\begin{lemma} Let $(A,\lambda,\rho,V)$ be a point on $X_U(k)$ together
with a $(\tau_0,i)$-special $V$, and let 
$(B,\lambda^{\prime},f \circ \rho \circ \phi^{-1})$ be the corresponding
point of $X_{U^{\prime}}(k)$.  Let 
$$W = \ker f: (\bCD_B)_{p_{\tau_0}} \rightarrow (\bCD_A)_{p_{\tau_0}}.$$
Then $W$ is $(\tau_0,i)$-constrained.
\end{lemma}
\begin{proof}
Note that since $f: (\bCD_B)_{p_{\tau}} \rightarrow (\bCD_A)_{p_{\tau}}$
is an isomorphism for $\tau \neq \tau_0$, we have that
$$W = \ker f: \bigoplus_{\tau} (\bCD_B)_{p_{\tau}} \rightarrow 
\bigoplus_{\tau} (\bCD_A)_{p_{\tau}}.$$
In particular $W$ is stable under $F$ and $V$; but since $F$ and $V$
send $W$ to $(\bCD_{p_{\sigma\tau_0}})$ and $(\bCD_{p_{\sigma^{-1}\tau_0}})$,
and neither of these contain any nonzero element of $W$, we have that
$W$ is killed by both $F$ and $V$.  The result follows immediately.
\end{proof}

We have thus associated to each tuple $(A,\lambda,\rho, V)$
a tuple $(B,\lambda^{\prime}, f\circ \rho \circ \phi^{-1}, W).$
We now describe an inverse construction.

Let $(B,\lambda^{\prime}, \rho^{\prime})$ be a point in $X_{U^{\prime}}(k)$,
and let $W$ be a $(\tau_0,i)$-constrained subspace of $(\bCD_B)_{p_{\tau_0}}$.
Define a submodule $N_W$ of $\bCD_B$ by:
\begin{enumerate}
\item $(N_W)_{p_{\tau_0}} = W$
\item $(N_W)_{\tau} = 0$ for $\tau \neq \tau_0$
\item $(N_W)_{q_{\tau}} = (N_W)_{p_{\tau}}^{\perp}$ for all $\tau$.
\end{enumerate}

It is clear that $N_W$ is stable under $F$ and $V$, and is a maximal
isotropic submodule of $\bCD_B$.  The inclusion of $N_W$ in $\bCD_B$
fits into an exact sequence
$$0 \rightarrow N_W \rightarrow \bCD_B \rightarrow \bCD(K^{\prime}) 
    \rightarrow 0,$$
where $\bCD(K^{\prime})$ is the Dieudonn{\'e} module of a subgroup
$K^{\prime}$ of $B$. 

Let $A = B/K^{\prime}$, and let $f^{\prime}: B \rightarrow A$ be the 
natural quotient map.  Then, just as before, there is a natural
polarization $\lambda$ on $A$ such that $p\lambda = (f^{\prime})^{\vee}
\lambda^{\prime} f^{\prime}$.

\begin{lemma} The dimension of $\Lie(A/k)_{p_{\tau}}$ 
(resp. $\Lie(A/k)_{q_{\tau}}$) is $r_{\tau}$ (resp. $s_{\tau}$) for
all $\tau$.   
\end{lemma}

\begin{proof}
The proof of this lemma is identical to the proof of Lemma~\ref{lemma:lierank},
and we omit it.
\end{proof}

It follows that the triple 
$(A,\lambda, \frac {1}{p}f^{\prime} \circ \rho^{\prime} \circ \phi)$
is a $k$-valued point of $X_U$.  Moreover, define $H$ by
$$H = \ker f^{\prime}: (\bCD_A)_{p_{\tau_0}} \rightarrow 
(\bCD_B)_{p_{\tau_0}}.$$  Then we have:

\begin{lemma} The space $H$ is $(\tau_0,i)$-special.
\end{lemma}
\begin{proof}
Since the image of $f^{\prime}: \bCD_A \rightarrow \bCD_B$ is
$N_W$, and $(N_W)_{p_{\tau_0}}$ has dimension $i$, $V$ has dimension
$n-i$.  The submodule $M_H = \ker f^{\prime}: \bCD_A \rightarrow \bCD_B$
is stable under $F$ and $V$, so in particular $F((M_H)_{p_{\sigma^{-1}\tau_0}})$
is contained in $(M_H)_{p_{\tau_0}}$.  But the former is all of
$F((\bCD)_A)_{p_{\sigma^{-1}\tau_0}}$, whereas the latter is just $H$.
In particular $H$ contains $F((\bCD_A)_{p_{\sigma^{-1}\tau_0}})$.
Similarly $H$ contains $V((\bCD_A)_{p_{\sigma\tau_0}}),$ so
$H$ is $(\tau_0,i)$-special.
\end{proof}

\begin{theorem} \label{thm:points}
The constructions above that associate to each $(A,\lambda,\rho,H)$
the corresponding $(B,\lambda^{\prime},\rho^{\prime},W)$ (and vice versa) 
are inverse to each other.
In particular there is a natural bijection between the space of
tuples $(A,\lambda,\rho,H)$ where $(A,\lambda,\rho) \in X_U(k)$
and $H$ is $(\tau_0,i)$-special, and the space of
tuples $(B,\lambda^{\prime},\rho^{\prime}, W)$ where
$(B,\lambda^{\prime},\rho^{\prime}) \in X_{U^{\prime}}(k)$ and
$W$ is $(\tau_0,i)$-constrained.
\end{theorem}
\begin{proof}
Fix a particular $(A,\lambda,\rho,H)$, and let $(B,\lambda^{\prime},
\rho^{\prime},W)$ be the point associated to it by the first construction
above.  Let $(A^{\prime\prime},\lambda^{\prime\prime},\rho^{\prime\prime},
H^{\prime\prime})$ be the point associated to 
$(B,\lambda^{\prime},\rho^{\prime},W)$ by the second construction above.

We need to show that the tuples $(A,\lambda,\rho,H)$ and
$(A^{\prime\prime},\lambda^{\prime\prime},\rho^{\prime\prime},H^{\prime\prime})$
are isomorphic.  Let $f: A \rightarrow B$ be the map used in the construction
of $B$ from $A$, and $f^{\prime}: B \rightarrow A^{\prime\prime}$ be
the map used in the construction of $A^{\prime\prime}$ from $B$.
The composition $f^{\prime}f$ induces the zero map $\bCD_{A^{\prime\prime}}
\rightarrow \bCD_A$, and hence its kernel contains $A[p]$.  Degree
considerations then show that the kernel is exactly $A[p]$, so
that $\frac{1}{p}f^{\prime}f$ is an isomorphism of $A$ with $A^{\prime\prime}$.
It is easy to check that this isomorphism carries $\lambda$ to
$\lambda^{\prime\prime}$ and $\rho$ to $\rho^{\prime\prime}$.  We
henceforth identify $A$ with $A^{\prime\prime}$ via this isomorphism.

Note now that by construction, we have
$$H^{\prime\prime} = \ker f^{\prime}: (\bCD_A)_{p_{\tau_0}} \rightarrow
(\bCD_B)_{p_{\tau_0}}.$$
By our definition of $f$, we have that $H = f((\bCD_B)_{p_{\tau_0}})$.
Since 
$$f(\bCD_B) = \ker f^{\prime}: \bCD_A \rightarrow \bCD_B,$$
it follows that $H = H^{\prime\prime}$.  Thus the second construction
is a left inverse to the first.

The proof that the second construction is a right inverse to the first
is similar, and will be omitted.
\end{proof}

\section{Geometrizing the Construction} \label{sec:geometry}

We now make our calculations with points in the previous section into
a geometric relationship between $X_{U^{\prime}}$ and $X_U$, by realizing
the bijection above as arising from a map of varieties.  We also study
the relationship of these varieties to $X_U$ and $X_{U^{\prime}}$.
We do so by systematically replacing the Dieudonn{\'e} modules appearing
in the previous section with DeRham cohomology modules.

\begin{definition} Let $S$ be a $k_0$-scheme, and 
$(A,\lambda, \rho)$ a point of $X_U(S)$.  A subbundle $H$
of $H^1_{DR}(A/S)_{p_{\tau_0}}$ is $(\tau_0,i)$-special if
$H$ has rank $n-i$, and contains both $\Lie(A/S)_{p_{\tau_0}}^*$ and 
$\Fr(H^1_{\DR}(A^{(p)}/S)_{p_{\tau_0}})$, where $\Fr$ denotes the relative
Frobenius $A \rightarrow A^{(p)}$.
\end{definition}

This generalizes our previous notion for the case when $S = \spec k,$
$k$ perfect.  

\begin{lemma} Let $(A,\lambda,\rho)$ be a point of $X_U(S)$.  Then 
$(A,\lambda,\rho)$ admits a $(\tau_0,i)$-special $H$ if and only if
the rank of 
$$\Ver: \Lie(A^{(p)}/S)_{p_{\tau_0}} \rightarrow \Lie(A/S)_{p_{\tau_0}}$$ 
is less than or equal to $r_{\sigma^{-1}\tau_0} - i$.
\end{lemma}
\begin{proof}
The kernel of 
$$\Ver: H^1_{DR}(A/S) \rightarrow H^1_{DR}(A^{(p)}/S)$$ 
is equal to the image of 
$$\Fr: H^1_{DR}(A^{(p)}/S) \rightarrow H^1_{DR}(A/S).$$
Since the dual of the map $\Ver: \Lie(A^{(p)}/S) \rightarrow \Lie(A/S)$
is the restriction of the map 
$\Ver: H^1_{DR}(A^{(p)}/S) \rightarrow H^1_{DR}(A/S)$
to the submodule $\Lie(A^{(p)}/S)^*$ of $H^1_{DR}(A^{(p)}/S)$, the rank of 
the map
$$\Ver: \Lie(A^{(p)}/S)_{p_{\tau_0}} \rightarrow \Lie(A/S)_{p_{\tau_0}}$$ 
is less than or equal to $r_{\sigma^{-1}\tau_0} - i$ if and only
if the rank of the intersection of the subsheaves
$\Lie(A/S)_{p_{\tau_0}}^*$ and $\Fr(H^1_{DR}(A^{(p)}/S)_{p_{\tau_0}})$
of $H^1_{DR}(A^{(p)})$ has rank at least $r_{\tau_0} - 
r_{\sigma^{-1}\tau_0} + i$.
This is true if and only if their sum has rank at most $n-i$, which in
turn is true if and only if there exists an subbundle $H$ of rank
$n-i$ containing both of them.
\end{proof}

Let $\CA$ denote the universal abelian variety on $X_U$.  We let
$(X_U)_{\tau_0,i}$ denote the subscheme of $X_U$ on which the map
$\Ver: \Lie(\CA^{(p)}/X_U)_{p_{\tau_0}} \rightarrow \Lie(\CA/X_U)_{p_{\tau_0}}$
has rank less than or equal to $r_{\sigma^{-1}\tau_0} - i$.  The
closed points on $(X_U)_{\tau_0,i}$ are precisely the $(\tau_0,i)$-special
points in the language of the preceding section.  (In particular, the
results of the previous section show that $(X_U)_{\tau_0,i}$ is nonempty.)

Let $\tX_{U,\tau_0,i}$ denote the $k_0$-scheme parametrizing
tuples $(A,\lambda,\rho,H)$, where $(A,\lambda,\rho) \in X_U(S)$
and $H$ is a $(\tau_0,i)$-special subspace of $H^1_{DR}(A/S)_{p_{\tau_0}}.$
There is a natural map $\tX_{U,\tau_0,i} \rightarrow X_U$, whose image is 
contained in $(X_U)_{\tau_0,i}$.

Our first goal is to understand the map $\tX_{U,\tau_0,i} \rightarrow X_U$.  
We will do so by constructing a local model for this map. 

For $\tau \neq \tau_0$, let $\CM_{\tau} = G(r_{\tau},n)_{\FF_p}$ be the 
Grassmannian parametrizing $r_{\tau}$-planes in $\FF_p^n$.  Define
$\CM_{\tau_0}$ to be the Schubert cycle in $G(r_{\tau_0},n)$ parametrizing
$r_{\tau_0}$-planes in $\FF_p^n$ that intersect the span of the first
$n-r_{\sigma^{-1}\tau}$ basis vectors in $\FF_p^n$ in a subspace of
dimension at least $r_{\tau_0} - r_{\sigma^{-1}\tau_0} + i$.

Finally, define $\tCM_{\tau_0}$ to be the moduli space parametrizing
pairs $(V,H)$, where $V$ is a subspace of $\FF_p^n$ of dimension
$r_{\tau_0}$, and $H$ is a subspace of $\FF_p^n$ of dimension $n-i$
containing both $V$ and the span of the first $n-r_{\sigma^{-1}\tau_0}$
basis vectors in $\FF_p^n$.  There is a natural map $\tCM_{\tau_0} 
\rightarrow \CM_{\tau_0}$ that forgets $H$; 
this map is generically one-to-one.

On the other hand, we have a natural map $\tCM_{\tau_0} \rightarrow
G(n-i,n)_{\FF_p}$ that forgets $V$.  The fibers of this map over a given
$H$ are simply $G(r_{\tau_0},H)$.  It follows that $\tCM_{\tau_0}$ is smooth;
it is a natural desingularization of $\CM_{\tau_0}$.

For $\tau \neq \tau_0$, set $\tCM_{\tau} = \CM_{\tau}$.
Let $\CM$ be the product (over $\FF_p$) of the $\CM_{\tau}$ for all $\tau$,
and let $\tCM$ be the product of the $\tCM_{\tau}$ for all $\tau$. 
We have a natural map $\tCM \rightarrow \CM$.

\begin{theorem} \label{thm:locmodel}
The map $\tCM \rightarrow \CM$ is a local model for the map
$\tX_{U,\tau_0,i} \rightarrow (X_U)_{\tau_0,i}$, in the sense that for 
any field $k$, and every $x \in (X_U)_{\tau_0,i}(k)$, there is a point 
$p$ of $\CM(k)$ and {\'e}tale neighborhoods $U_x$ of $x$ and $U_p$ of $p$ such 
that the base change of $\tX_U \rightarrow (X_U)_{\tau_0,i}$ to $U_x$ is 
isomorphic to the base change of $\tCM \rightarrow \CM$ to $U_p$.
\end{theorem}

To prove this, we first introduce two new schemes $(X_U)_{\tau_0,i}^+$ and
$\tX^+_{U,\tau_0,i}$.  The former parametrizes tuples 
$(A,\lambda,\rho,\{e_{i,\tau}\})$, where $i$ runs from $1$ to $n$ for each
$\tau: F^+ \rightarrow \RR$, and the set $\{e_{1,\tau}, \dots, e_{n,\tau}\}$
is a basis for $H^1_{\DR}(A)_{p_{\tau}}$ for all $\tau$, such that
the subset $\{e_{1,\tau_0}, \dots, e_{n-r_{\sigma^{-1}\tau_0},\tau_0}\}$
of $\{e_{1,\tau_0}, \dots, e_{n,\tau_0}\}$
is a basis for the subbundle $\Fr(H^1_{\DR}(A^{(p)})_{p_{\tau_0}})$
of $H^1_{\DR}(A)_{p_{\tau_0}}$.  The latter parametrizes the same data,
plus a $(\tau_0,i)$-special subbundle $H$ of $H^1_{\DR}(A)_{p_{\tau_0}}$.

Clearly $(X_U)_{\tau_0,i}^+$ and $\tX^+_{U,\tau_0,i}$ possess natural maps
to $(X_U)_{\tau_0,i}$ and $\tX_{U,\tau_0,i}$, respectively, by forgetting the
$e_{i,\tau}$.  They also possess natural maps to $\CM$ and $\tCM$, which
we will now construct. 

Given an $S$-valued point $(A,\lambda,\rho,\{e_{i,\tau}\})$ of
$(X_U)_{\tau_0,i}$, the basis $e_{i,\tau}$ allows us to identify
$H^1_{\DR}(A)_{p_{\tau}}$ with $\OO_S^n$.  Then the subbundle
$\Lie(A/S)_{p_{\tau}}^*$ of $H^1_{\DR}(A)_{p_{\tau}}$ gives us a
corresponding subbundle $V$ of $\OO_S^n$, and hence a point of $\CM_{\tau}$.
We thus obtain a map $(X_U)^+_{\tau_0,i} \rightarrow \CM$.  If in
addition we have a $(\tau_0,i)$-special subbundle $H$ of
$H^1_{\DR}(A)_{p_{\tau_0}}$, then the pair $(\Lie(A)^*,H)$ corresponds
to a point of $\tCM_{\tau_0}$.  We therefore obtain a map 
$\tX^+_{U,\tau_0,i} \rightarrow \tCM$.  These fit into a commutative diagram:
$$
\begin{array}{ccccc}
\tX_{U,\tau_0,i} & \leftarrow & \tX^+_{U,\tau_0,i} & \rightarrow & \tCM \\
\downarrow & & \downarrow & & \downarrow \\
(X_U)_{\tau_0,i} & \leftarrow & (X_U)^+_{\tau_0,i} & \rightarrow & \CM.
\end{array}
$$

The left-hand horizontal maps are clearly smooth.  We will show the right-hand
horizontal maps are also smooth.

The right-hand
square in the above diagram is cartesian, so it suffices to show that
the map $(X_U)^+_{\tau_0,i} \rightarrow \CM$ is smooth.  There is a standard
way to do this
using the crystalline deformation theory of abelian varieties.  We first
summarize the necessary facts:

Let $S$ be a scheme, and $S^{\prime}$
a thickening of $S$ equipped with divided powers.  Let $\CC_{S^{\prime}}$
denote the category of abelian varieties over $S^{\prime}$, and
$\CC_S$ denote the category of abelian varieties over $S$.  For $A$ an
object of $\CC_{S^{\prime}}$, let $\overline{A}$ denote its base change
to $\CC_S$.

Fix an $A$ in $\CC_{S^{\prime}}$, and consider the module
$H^1_{\cris}(\overline{A}/S)_{S^{\prime}}$.  This is a locally free
$\OO_{S^{\prime}}$-module, and we have a canonical isomorphism:
$$H^1_{\cris}(\overline{A}/S)_{S^{\prime}} \cong H^1_{\DR}(A/S^{\prime}).$$
Moreover, we have a natural submodule
$$\Lie(A/S^{\prime})^* \subset H^1_{\DR}(A/S^{\prime}).$$
The preceding isomorphism thus gives us a subbundle of
$H^1_{\cris}(\overline{A}/S)_{S^{\prime}}$ that lifts the subbundle 
$\Lie(\overline{A}/S)^*$ of $H^1_{\DR}(\overline{A}/S)$.

Knowing this lift allows us to recover $A$ from $\overline{A}$.
More precisely, let $\CC_S^+$ denote the category of pairs
$(\overline{A},\omega)$, where $\overline{A}$ is an object of $\CC_S$
and $\omega$ is a subbundle of 
$H^1_{\cris}(\overline{A}/S)_{S^{\prime}}$ that lifts
$\Lie(\overline{A}/S)^*$.  Then the construction outlined above gives us
a functor from $\CC_{S^{\prime}}$ to $\CC_S^+$.

\begin{theorem}[Grothendieck] \label{thm:grot}
The functor $\CC_{S^{\prime}} \rightarrow \CC_S^+$ defined above is an
equivalence of categories.
\end{theorem}

\begin{proof} The proof is sketched in~\cite{montreal}, pp. 116-118.
A complete proof can be found in~\cite{mazurmessing}.
\end{proof}

\begin{proposition}
The map $(X_U)^+_{\tau_0,i} \rightarrow \CM$ is smooth.
\end{proposition}
\begin{proof}
Let $R^{\prime}$ be
a ring, and $I$ an ideal of $R$ such that $I^2 = \{0\}$.  Let $R$ be the
ring $R^{\prime}/I$.  It suffices to show for any diagram
$$
\begin{array}{ccc}
\spec R & \rightarrow & (X_U)^+_{\tau_0,i} \\
\downarrow & & \downarrow \\
\spec R^{\prime} & \rightarrow & \CM
\end{array}$$
there is a map $\spec R^{\prime} \rightarrow (X_U)^+_{\tau_0,i}$.

In terms of the moduli, such a diagram consists of the
following data:
\begin{enumerate}
\item an $R$-valued point $(A,\lambda,\rho)$ of $(X_U)^+_{\tau_0,i}$,
\item for each $\tau$, bases $e_{i,\tau}$ of $H^1_{\DR}(A/R)_{p_{\tau}}$,
such that the set $e_{1,\tau_0}, \dots, e_{s_{\sigma^{-1}\tau,\tau_0}}$
is a basis for the submodule $\Fr(H^1_{\DR}(A^{(p)}/R)_{p_{\tau_0}})$ of
$H^1_{\DR}(A/R)_{p_{\tau_0}}$,
\item For each $\tau$, a rank $r_{\tau}$ subbundle $V_{\tau}$
of $(R^{\prime})^n$ whose reduction modulo $I$ is the subbundle of
$R^n$ that corresponds to the subbundle $\Lie(A/R)_{p_{\tau}}^*$ of
$H^1_{\DR}(A/R)_{p_{\tau}}$ under the identification of the latter with
$R^n$ induced by the $e_{i,\tau}$.  The bundle $V_{\tau_0}$ has the
additional property that its intersection with the span of the
first $s_{\sigma^{-1}\tau}$ standard basis vectors of $(R^{\prime})^n$
has rank at least $r_{\tau_0} - r_{\sigma^{-1}\tau_0} + i$.
\end{enumerate}

For each $\tau$ and $i$, let $\te_{\tau,i}$ be a lift of $e_{\tau,i}$
to $(H^1_{\cris}(A/R)_{R^{\prime}})_{p_{\tau}}$.  (If $\tau = \tau_0$
and $i \leq s_{\sigma^{-1}\tau_0}$, then we require that $\te_{\tau_i}$ lies
in the subbundle $\Fr(H^1_{\cris}(A^{(p)}/R)_{R^{\prime}})_{p_{\tau}}$
of $(H^1_{\cris}(A/R)_{R^{\prime}})_{p_{\tau}}$.)

Under this choice of basis, each $V_{\tau}$ corresponds to a subbundle
$\omega_{p_{\tau}}$ of $(H^1_{\cris}(A/R)_{R^{\prime}})_{p_{\tau}}$ that
lifts the subbundle $\Lie(A/R)^*_{p_{\tau}}$ of $H^1_{\DR}(A/R)_{p_{\tau}}$.
Define $\omega_{q_{\tau}} = \omega_{p_{\tau}}^{\perp}$ for all $\tau$,
where $\perp$ denotes orthogonal complement with respect to the pairing
$$(H^1_{\cris}(A/R)_{R^{\prime}})_{p_{\tau}} \times 
(H^1_{\cris}(A/R)_{R^{\prime}})_{q_{\tau}} \rightarrow R^{\prime}$$
induced by $\lambda$.

By crystalline deformation theory, this defines a lift of $A$ to an abelian
scheme over $\spec R^{\prime}$.  The relation $\omega_{q_{\tau}}=
\omega_{p_{\tau}}^{\perp}$ implies that $\lambda$ lifts to a 
prime-to-$p$ polarization of this lift as well.  We thus obtain
a point $(\tA,\tlambda,\trho)$ of $X_U(R^{\prime})$.  Moreover,
since the rank of the intersection of $V_{\tau_0}$ with the span of
the first $s_{\sigma^{-1}\tau_0}$ basis vectors of $(R^{\prime})^n$
has rank at least $r_{\tau} - r_{\sigma^{-1}\tau} + i$, the
same can be said for the intersection of $\omega_{p_{\tau_0}}$
with $\Fr(H^1_{\cris}(A/R)_{R^{\prime}})_{p_{\tau_0}}$, and
hence also for the intersection of $\Lie(\tA/R^{\prime})^*_{p_{\tau_0}}$
with $\Fr(H^1_{\DR}(A/R)_{p_{\tau_0}})$.  Thus $(\tA,\tlambda,\trho)$
lies in $(X_U)_{\tau_0,i}$.  Finally, the basis $\te_{\tau_i}$
corresponds to a basis of $H^1_{\DR}(\tA/R^{\prime})_{p_{\tau}}$ for
each $\tau$, and these bases, together with the point
$(\tA,\tlambda,\trho)$ define the required point of $(X_U)^+_{\tau_0,i}$. 
\end{proof}

It is easy to see (for instance, by computing the dimension of the
tangent space to a fiber) that the smooth maps
$(X_U)^+_{\tau_0,i} \rightarrow (X_U)_{\tau_0,i}$ and
$(X_U)^+_{\tau_0,i} \rightarrow \CM$ have the same relative dimension.
Thus if $x$ is a point of $(X_U)_{\tau_0,i}$, $x^+$ is a lift of $x$
to $(X_U)^+_{\tau_0,i}$, and $p$ is the image of $x^+$ in $\CM$, the
complete local ring ${\hat \OO}_{(X_U)^+_{\tau_0,i}, x^+}$ is simultaneously
a power series ring over ${\hat \OO}_{(X_U)_{\tau_0,i}, x}$ and a power
series ring over ${\hat \OO}_{\CM,p}$, in the same number of variables.

Corollary 4.6 of~\cite{sgee} then implies that 
${\hat \OO}_{(X_U)_{\tau_0,i},x}$ and ${\hat \OO}_{\CM,p}$ are isomorphic.
More precisely, the proof of this corollary shows that there is a
map ${\hat \OO}_{\CM,p} \rightarrow {\hat \OO}_{(X_U)^+_{\tau_0,i},x^+}$
whose composition with the map ${\hat \OO}_{(X_U)^+_{\tau_0,i},x^+}
\rightarrow {\hat \OO}_{(X_U)_{\tau_0,i},x}$ is an isomorphism, and
whose composition with the map ${\hat \OO}_{(X_U)^+_{\tau_0,i},x^+}
\rightarrow {\hat \OO}_{\CM,p}$ is the identity on ${\hat \OO}_{\CM,p}.$

It follows by Artin approximation (\cite{artin}, especially Corollary 2.5) 
that there are {\'e}tale neighborhoods $U_x$, $U_{x^+}$, and $U_p$
of $x$, $x^+$, and $p$ respectively, a diagram
$$
\begin{array}{ccccc}
U_x & \leftarrow & U_{x^+} & \rightarrow & U_p\\
\downarrow & & \downarrow & & \downarrow\\
(X_U)_{\tau_0,i} & \leftarrow & (X_U)^+_{\tau_0,i} & \rightarrow &
\CM
\end{array}
$$
in which both squares are Cartesian, 
and a section $U_p \rightarrow U_{x^+}$ whose composition
with the map $U_{x^+} \rightarrow U_x$ is an isomorphism, and whose
composition with the map $U_{x^+} \rightarrow U_p$ is the identity on $U_p$.

Define 
$$\tU_x = \tX_{U,\tau_0,i} \times_{(X_U)_{\tau_0,i}} U_x,$$
$$\tU_{x^+} = \tX^+_{U,\tau_0,i} \times_{(X_U)^+_{\tau_0,i}} U_{x^+},$$ 
$$\tU_p = \tCM \times_{\CM} U_p.$$  
We obtain from the section $U_p \rightarrow U_{x^+}$
a map $\tU_p \rightarrow \tU_{x^+}$ whose composition with the natural map
$\tU_{x^+} \rightarrow \tU_x$ is an isomorphism.  This yields a commutative
square
$$ 
\begin{array}{ccc}
\tU_x & \cong & \tU_p\\
\downarrow & & \downarrow\\
U_x & \cong & U_p,
\end{array}
$$
and thus establishes Theorem~\ref{thm:locmodel}.

Theorem~\ref{thm:locmodel} implies that the singularities of
$(X_U)_{\tau_0,i}$ look ({\'e}tale locally) like products of an
affine space with a singularity of the Schubert cycle $\CM_{\tau_0}$.
Moreover, $\tX_{U,\tau_0,i}$ is a natural desingularization
of $(X_U)_{\tau_0,i}$.  For $j \geq 0$, the fiber of the map
$\tX_{U,\tau_0,i} \rightarrow (X_U)_{\tau_0,i}$ over a point of
$(X_U)_{\tau_0,i+j} \setminus (X_U)_{\tau_0,i+j+1}$ is a Grassmannian
parametrizing $j$-planes in an $i+j$-dimensional space.

The points of $\tX_{U,\tau_0,i}$ over a perfect field $k$ correspond to
tuples $(A,\lambda,\rho,H)$, where $(A,\lambda,\rho)$ is a $k$-valued
point of $X_U$, and $H$ is a $(\tau_0,i)$-special subspace of 
$\CD(A[p])_{p_{\tau_0}}$.  In order to geometrize the construction
in the previous section, we would like to have a map from
$\tX_{U,\tau_0,i}$ to $X_{U^{\prime}}$.
Unfortunately, $\tX_{U,\tau_0,i}$ does not admit such a map.
We must therefore introduce another moduli problem:

\begin{definition}
Let $S$ be a $k_0$-scheme, $(A,\lambda,\rho)$ a point
of $X_U(S)$, and $(B,\lambda^{\prime},\rho^{\prime})$ a point of
$X_{U^{\prime}}(S)$.  A $(\tau_0,i)$-special isogeny $f: (A,\lambda,\rho)
\rightarrow (B,\lambda^{\prime},\rho^{\prime})$ is an $\OO_F$-isogeny $f: A
\rightarrow B$, of degree $p^{nd}$, such that:
\begin{enumerate}
\item $p\lambda = f^{\vee}\lambda^{\prime} f$,
\item the $U^{\prime}$-level structure $\rho^{\prime}$ on $B$ corresponds to
$f \circ \rho$ under the identification of $T$ with $T^{\prime}$
fixed in the previous section,
\item for each $\tau \neq \tau_0$, the map $f$ induces an isomorphism
of $H^1_{\DR}(B/S)_{p_{\tau}}$ with $H^1_{\DR}(A/S)_{p_{\tau}}$, and
\item the image of $H^1_{\DR}(B/S)_{p_{\tau_0}}$ in 
$H^1_{\DR}(A/S)_{p_{\tau_0}}$ under $f$ has rank $n-i$.  (It is necessarily
a subbundle of $H^1_{\DR}(A/S)_{p_{\tau_0}}$.)
\end{enumerate}
We denote by $\hX_{U,\tau_0,i}$ the scheme parametrizing tuples
$(A,\lambda,\rho,B,\lambda^{\prime},\rho^{\prime},f)$, where
$(A,\lambda,\rho)$ is a point of $(X_U)$, $(B,\lambda^{\prime},\rho^{\prime})$
is a point of $X_{U^{\prime}}$, and $f$ is a $(\tau_0,i)$-special isogeny
from $(A,\lambda,\rho)$ to $(B,\lambda^{\prime},\rho^{\prime})$.
\end{definition}

If $(A,\lambda,\rho,B,\lambda^{\prime},\rho^{\prime},F)$ is
a point of $\hX_{U,\tau_0,i}(S)$, then $f(H^1_{\DR}(B/S)_{p_{\tau_0}})$
is a $(\tau_0,i)$-special subbundle of $H^1_{\DR}(A/S)_{p_{\tau_0}}$.
Indeed, we know that the kernel of
$$f: H^1_{\DR}(B/S)_{p_{\tau_0}} \rightarrow H^1_{\DR}(A/S)_{p_{\tau_0}}$$
has rank $i$.  The subbundle $\Lie(B/S)^*_{p_{\tau_0}}$ of 
$H^1_{\DR}(B/S)_{p_{\tau_0}}$ has rank $r_{\tau_0} + i$, and 
$\Lie(A/S)_{p_{\tau_0}}$ has rank $r_{\tau_0}$.  As $f$ maps the former
to the latter, $f(H^1_{\DR}(B/S)_{p_{\tau_0}})$ must contain 
$\Lie(A/S)^*_{p_{\tau_0}}.$  An identical argument shows that 
$f(H^1_{\DR}(B/S)_{p_{\tau_0}})$ contains 
$\Fr(H^1_{\DR}(A^{(p)}/S)_{p_{\tau_0}})$.  The morphism of functors that
associates the tuple $(A,\lambda,\rho,f(H^1_{\DR}(B/S)_{p_{\tau_0}}))$ to
the tuple $(A,\lambda,\rho,B,\lambda^{\prime},\rho^{\prime},f)$ therefore
induces a map $\hX_{U,\tau_0,i} \rightarrow \tX_{U,\tau_0,i}$.

\begin{proposition} The map $\hX_{U,\tau_0,i} \rightarrow \tX_{U,\tau_0,i}$ 
is a bijection on $k$-valued points for any perfect field $k$.
\end{proposition}
\begin{proof}
The construction in the previous section associates to every
$(A,\lambda,\rho)$ in $X_U(k)$, and every $(\tau_0,i)$-special
subspace $H$ of $\CD(A[p])_{p_{\tau_0}}$ (or equivalently of
$H^1_{\DR}(A/k)_{p_{\tau_0}}$) a $(B,\lambda^{\prime},\rho^{\prime})$
and a $(\tau_0,i)$-special isogeny 
$f: (A,\lambda,\rho) \rightarrow (B,\lambda^{\prime},\rho^{\prime}).$
This construction is inverse to the map 
$$\hX_{U,\tau_0,i}(k) \rightarrow \tX_{U,\tau_0,i}(k)$$
constructed above.
\end{proof}

This has strong consequences for the geometry of the map
$\hX_{U,\tau_0,i} \rightarrow \tX_{U,\tau_0,i}.$  In particular
we have the following result:

\begin{proposition} \label{prop:insep}
Let $Y$ and $Z$ be schemes of finite type over a perfect
field $k$ of characteristic $p$, such that $Z$ is normal and $Y$ is reduced.
Let
$f: Y \rightarrow Z$ be a proper map that is a bijection on points.
Then there is a map 
$f^{\prime}: Z_{p^r} \rightarrow Y$ such that
$$ff^{\prime}: Z_{p^r} \rightarrow Z$$ is the $r$th
power of the Frobenius.  (In particular $f$ is an isomorphism on {\'e}tale 
cohomology.)
\end{proposition}

This is proven in~\cite{u2}, Proposition 4.8.

\begin{remark} \rm One might wonder if the map $\hX_{U,\tau_0,i}
\rightarrow \tX_{U,\tau_0,i}$ is actually an isomorphism,
but in fact a straightforward calculation, using
Theorem~\ref{thm:grot}, shows that it often fails to be
an isomorphism on tangent spaces.
\end{remark}

The scheme $\hX_{U,\tau_0,i}$ admits an obvious map to $X_{U^{\prime}}$.
In fact, as one expects from the previous section, it admits
a map to a scheme $(X_{U^{\prime}})^{\tau_0,i}$ parametrizing
$(\tau_0,i)$-constrained subspaces.  More precisely:

\begin{definition}
Let $S$ be a $k_0$-scheme, and
$(B,\lambda^{\prime},\rho^{\prime})$ a point of $X_{U^{\prime}}(S)$.
A subbundle $W$ of $H^1_{\DR}(B/S)_{p_{\tau_0}}$ is $(\tau_0,i)$-constrained
if it has rank $i$ and is contained in both $\Lie(B/S)^*_{p_{\tau_0}}$
and $\Fr(H^1_{\DR}(B^{(p)}/S)_{p_{\tau_0}})$.  We denote by
$(X_{U^{\prime}})^{\tau_0,i}$ the scheme parametrizing tuples
$(B,\lambda^{\prime},\rho^{\prime},W)$, where 
$(B,\lambda^{\prime},\rho^{\prime})$ is a point of $X_{U^{\prime}}$
and $W$ is a $(\tau_0,i)$-constrained subbundle of 
$H^1_{\DR}(B/S)_{p_{\tau_0}}.$
\end{definition}

\begin{proposition}
Let $(A,\lambda,\rho,B,\lambda^{\prime},\rho^{\prime},f)$ be an
element of $\hX_{U,\tau_0,i}(S)$, and let $W$ be the kernel of
the map 
$$f: H^1_{\DR}(B/S)_{p_{\tau_0}} \rightarrow H^1_{\DR}(A/S)_{p_{\tau_0}}.$$
Then $W$ is a $(\tau_0,i)$-constrained subbundle of 
$H^1_{\DR}(B/S)_{p_{\tau_0}}.$
\end{proposition}
\begin{proof}
The rank of $W$ is clearly $i$, so it suffices to show that $W$ is contained in
$\Lie(B/S)^*_{p_{\tau_0}}$ and $\Fr(H^1_{\DR}(B^{(p)}/S)_{p_{\tau_0}})$.
The former has rank $r_{\tau_0} + i$, whereas $\Lie(A/S)^*_{p_{\tau_0}}$
has rank $r_{\tau_0}$.  Thus the kernel of the map
$$f:\Lie(B/S)^*_{p_{\tau_0}} \rightarrow \Lie(A/S)^*_{p_{\tau_0}}$$ 
has dimension at least $i$.  Since this kernel is contained in $W$, it
must be equal to $W$, and hence $W$ is contained in $\Lie(B/S)^*_{p_{\tau_0}}$.
The proof of containment in $\Fr(H^1_{\DR}(B/S)_{p_{\tau_0}})$ is
similar.
\end{proof}

We thus have a map $\hX_{U,\tau_0,i} \rightarrow (X_{U^{\prime}})^{\tau_0,i}$
that takes $(A,\lambda,\rho,B,\lambda^{\prime},\rho^{\prime},f)$
to $(B,\lambda^{\prime},\rho^{\prime},W)$, with $W$ as above.  For
any perfect field $k$ of characteristic $p$, composing the map
$$\hX_{U,\tau_0,i}(k) \rightarrow (X_{U^{\prime}})^{\tau_0,i}(k)$$
with the bijection
$$\tX_{U,\tau_0,i}(k) \rightarrow \hX_{U,\tau_0,i}(k)$$
yields the bijection
$$\tX_{U,\tau_0,i}(k) \rightarrow (X_{U^{\prime}})^{\tau_0,i}(k)$$
constructed in the previous section.  In particular the map
$\hX_{U,\tau_0,i} \rightarrow (X_{U^{\prime}})^{\tau_0,i}$ is
a bijection on points. 

\begin{lemma} The scheme $(X_{U^{\prime}})^{\tau_0,i}$ is smooth
over $k_0$.
\end{lemma}
\begin{proof}
The dimension of $(X_{U^{\prime}})^{\tau_0,i}$ is equal
to that of $\tX_{U,\tau_0,i}$, and hence to that of $\CM$.  Thus
$(X_{U^{\prime}})^{\tau_0,i}$ has dimension equal to
$$(\sum_{\tau} r_{\tau}s_{\tau}) - i(i + r_{\tau_0} - r_{\sigma^{-1}\tau_0}).$$

We must show that the dimension of the tangent space to
$(X_{U^{\prime}})^{\tau_0,i}$ at any $k$-valued point $x$ is equal
to this number.  Let $(B,\lambda^{\prime},\rho^{\prime},W)$ be the 
moduli object
corresponding to $x$, and let $S = \spec k[\epsilon]/\epsilon^2$.
Then, by Grothendieck's theorem, specifying a tangent vector to
$(X_{U^{\prime}})^{\tau_0,i}$ at $x$ is equivalent to specifying
the following data:
\begin{enumerate}
\item For each $\tau$, a lift $\omega_{p_{\tau}}$ of $\Lie(B/k)^*_{p_{\tau}}$
from $H^1_{\DR}(B/k)_{p_{\tau}}$ to $(H^1_{\cris}(B/k)_S)_{p_{\tau}}$, and
\item a lift $\tW$ of $W$ to a subspace of $(H^1_{\cris}(B/k)_S)_{p_{\tau_0}}$
that is contained in $\omega_{p_{\tau_0}}$ and in
$\Fr(H^1_{\cris}(B^{(p)}/k)_S)_{p_{\tau_0}}.$
\end{enumerate}

The space of possible lifts of $W$ that are contained in
$\Fr(H^1_{\cris}(B^{(p)}/k)_S)_{p_{\tau_0}}$ has dimension 
$is_{\sigma^{-1}\tau_0}$.
(Recall that $\Lie(B/k)^*_{p_{\sigma^{-1}\tau}}$ has dimension 
$r_{\sigma^{-1}\tau_0} - i$, so that $\Fr(H^1_{\DR}(B^{(p))}/k)_{p_{\tau_0}}$ 
and $\Fr(H^1_{\cris}(B^{(p)}/k)_S)_{p_{\tau_0}}$ have dimension 
$s_{\sigma^{-1}\tau_0} + i$.)  Once we have fixed such a lift, the space
of $\omega_{p_{\tau_0}}$ containing that lift has dimension 
$r_{\tau_0}(s_{\tau_0} - i)$, as $\Lie(B/k)^*_{p_{\tau_0}}$ has
dimension $r_{\tau_0} + i$.

On the other hand, $\Lie(B/k)^*_{p_{\sigma^{-1}\tau_0}}$ has
dimension $r_{\sigma^{-1}\tau_0} - i$, so the space of possible
$\omega_{p_{\sigma^{-1}\tau_0}}$ has dimension 
$(r_{\sigma^{-1}\tau_0} - i)(s_{\sigma^{-1}\tau_0} + i)$.  For $\tau$ not
equal to either $\tau_0$ or $\sigma^{-1}\tau_0$, the space of
possible $\omega_{p_{\tau}}$ has dimension $r_{\tau}s_{\tau}$.

Summing these, we find that the tangent space at $x$ has dimension
$$(\sum_{\tau} r_{\tau}s_{\tau}) - i(i + r_{\tau_0} - r_{\sigma^{-1}\tau_0}),$$
as desired.
\end{proof}

\begin{corollary}
The map
$\hX_{U,\tau_0,i}^{\red} \rightarrow (X_{U^{\prime}})^{\tau_0,i}$
induces an isomorphism on {\'e}tale cohomology.  
\end{corollary}
\begin{proof}
This is immediate from Proposition~\ref{prop:insep}.
\end{proof}

In summary, we have constructed a cycle $(X_U)_{\tau_0,i}$ on $X_U$,
and a natural desingularization $\tX_{U,\tau_0,i}$.  The geometry of
this desingularization is closely related to that of $X_{U^{\prime}}$;
in particular there is a scheme $(X_{U^{\prime}})^{\tau_0,i}$ defined
in terms of the universal abelian variety on $X_{U^{\prime}}$, that 
is ``nearly isomorphic'' to $\tX_{U,\tau_0,i}$, in the sense that
there exists a scheme $\hX_{U,\tau_0,i}^{\red}$ and a diagram:
$$\tX_{U,\tau_0,i} \leftarrow \hX_{U,\tau_0,i}^{\red} \rightarrow 
(X_{U^{\prime}})^{\tau_0,i}$$
in which both maps are bijections on points and isomorphisms on {\'e}tale 
cohomology.  In particular the {\'e}tale cohomology groups of 
$\tX_{U,\tau_0,i}$ and $(X_{U^{\prime}})^{\tau_0,i}$ are naturally 
isomorphic via these maps.

\section{Cohomology} \label{sec:cohomology}

We now explore the implications of the previous section for the {\'e}tale
cohomology of Shimura varieties.  Let $N = \sum_{\tau} r_{\tau}s_{\tau}$
be the dimension of $X_U$, and $r = i(i + r_{\tau_0} - r_{\sigma^{-1}\tau_0})$
be the codimension of $(X_U)_{\tau_0,i}$.  Then $X_{U^{\prime}}$ has
dimension $N-2r$.  For the purposes of this section, we consider
$X_U$, $X_{U^{\prime}}$, etc. as schemes over $\overline{\FF}_p$.

Fix an $\ell$ different from $p$, and let $\xi$ be a finite dimensional
$\overline{\QQ}_{\ell}$-representation of $G(\AA^f_{\QQ})$.  As 
in~\cite{HT}, III.2,
this determines a lisse $\overline{\QQ}_{\ell}$-sheaf $\CL_{\xi}$
on $X_U$.  Let $\xi^{\prime}$ be the representation of 
$G^{\prime}(\AA^f_{\QQ})$ induced by $\xi$ and our fixed isomorphism
of $G^{\prime}(\AA^f_{\QQ})$ with $G(\AA^f_{\QQ})$.  Then we also have a
lisse sheaf $\CL_{\xi^{\prime}}$ on $X_{U^{\prime}}$.  If
$\hat \pi$ (resp. ${\hat \pi}^{\prime}$) denotes the map
$\hX_{U,\tau_0,i}^{\red} \rightarrow X_U$ (resp. the map
$\hX_{U,\tau_0,i}^{\red} \rightarrow X_{U^{\prime}}$), then there is
a natural isomorphism ${\hat \pi}^*\CL_{\xi} \cong 
({\hat \pi}^{\prime})^*\CL_{\xi^{\prime}}$.
 
We will construct, for each $j$, a map 
$$H^j_{\et}(X_{U^{\prime}}, \CL_{\xi^{\prime}}) \rightarrow
H^{j+2r}_{\et}(X_U, \CL_{\xi}(r)).$$
Note that this takes the middle degree cohomology of $X_{U^{\prime}}$ to
the middle degree cohomology of $X_U$.  

In order to construct this map, let us consider the following situation.
Let $X$ and $Y$ be smooth over $\overline{\FF}_p$, of dimensions $N$
and $N-r$, respectively.  Let $\CF$ be a lisse sheaf on $X$, and
$\pi: Y \rightarrow X$ a proper map of $\overline{\FF}_p$-schemes.

Let $\theta_X$ and $\theta_Y$ denote the structure maps
$$\theta_X: X \rightarrow \spec \overline{\FF}_p$$
$$\theta_Y: Y \rightarrow \spec \overline{\FF}_p.$$
We have natural isomorphisms: 
$$R\theta_X^! \QQ_{\ell} \cong \QQ_{\ell}[2N](N)$$ 
$$R\theta_Y^! \QQ_{\ell} \cong \QQ_{\ell}[2(N-r)](N-r).$$
Since $R\theta_Y^! = R\pi^!R\theta_X^!$,
it follows that $R\pi^!\QQ_{\ell} \cong \QQ_{\ell}[-2r](-r)$.

It follows that for any lisse sheaf $\CF$ on $X$, $R\pi^!\CF$
is naturally isomorphic to $\pi^*\CF[-2r](-r)$; this isomorphism is
simply the tensor product of the above isomorphism with $\CF$.

Since $\pi$ is proper, we have $\pi_! = \pi_*$.  We thus have a unit
map
$$R\pi_*R\pi^!\CF \rightarrow \CF$$
and therefore a morphism
$$R\pi_*\pi^*\CF \rightarrow \CF[2r](r)$$
in the derived category of sheaves on $X$.  This induces a map
$$\nu_{\pi}: H^j_{\et}(Y,\pi^*\CF) \cong H^j_{\et}(X,R\pi_*\pi^*\CF)
\rightarrow H^{j+2r}_{\et}(X,\CF(r)).$$

Note that if $X$ and $Y$ are proper, $\nu_{\pi}$ is simply the Poincar{\'e}
dual of the natural map $H^{2n-2r-j}_{\et}(X,\CF) \rightarrow
H^{2n-2r-j}_{\et}(Y,\pi^*\CF)$, suitably Tate twisted.  On the other
hand, if $\pi$ is a closed immersion, $\nu_{\pi}$ is simply the Gysin
map.

This construction is compatible with cycle maps in the following
sense: we have a commutative diagram
$$
\begin{array}{ccc}
A^j(Y) & \rightarrow & A^{j+r}(X)\\
\downarrow & & \downarrow\\
H^{2j}(Y,\QQ_{\ell}(2j)) & \rightarrow & H^{2j+2r}(X,\QQ_{\ell}(2j+2r))
\end{array}
$$
where the vertical maps are cycle class maps and the map
$A^j(Y) \rightarrow A^{j+r}(X)$ is the proper pushforward $\pi_*$ of
cycles on $Y$ to cycles on $X$.

Although the map $\nu_{\pi}$ is difficult to describe directly,
we do have the following result:

\begin{lemma} \label{lemma:reflect}
Let $c_{\pi}$ denote the class of $H^{2r}_{\et}(Y,\QQ_{\ell}(r))$
associated to the cycle class $\pi^!\pi_*[Y] \in A^r(Y)$, where
$[Y]$ denotes the fundamental class of $Y$ in $A^0(Y)$ and
$\pi^!$ denotes the refined Gysin homomorphism $A^r(X) \rightarrow A^r(Y)$
of~\cite{Fulton}, Definition 8.1.2.  Let $\eta_{\CF}$ be
the composition of $\nu_{\pi}$ with the natural map
$$H^{j+2r}(X,\CF(r)) \rightarrow H^{j+2r}(Y,\pi^*\CF(r))$$
Then $\eta_{\CF}$ is given by cup product with the cohomology class $c_{\pi}$.
\end{lemma}
\begin{proof}
The map $\eta_{\CF}$ is induced by the morphism 
(in the derived category)
$$\xi_{\CF}: \pi^*\CF \rightarrow \pi^*\CF[-2r](r)$$
that is the composition of the sequence of morphisms:
$$
\begin{array}{rcl}
\pi^*\CF & \rightarrow & \pi^*R\pi_*\pi^*\CF \\
         & \rightarrow & \pi^*R\pi_*\pi^!\CF[2r](r) \\
         & \rightarrow & \pi^*\CF[2r](r).
\end{array}
$$
If we identify $\pi^*\CF$ with $\pi^*\CF \otimes \QQ_{\ell}$,
the map $\xi_{\CF}$ is simply $\id \otimes \xi_{\QQ_{\ell}}$.  It
follows that if $a \in H^j_{\et}(Y,\pi^*\CF)$, and 
$b \in H^{j^{\prime}}(Y,\QQ_{\ell})$, then $\eta_{\CF}(a \cup b) =
a \cup \eta_{\QQ_{\ell}}(b)$.  Taking $b = 1$ we see that
$\eta_{\CF}(a) = a \cup \eta_{\QQ_{\ell}}(1)$.  It thus remains to
compute $\eta_{\QQ_{\ell}}(1)$.

We have a commutative diagram:
$$
\begin{array}{ccccc}
A^0(Y) & \stackrel{\pi_*}{\rightarrow} & A^r(X) & \stackrel{\pi^!}{\rightarrow}
& A^r(Y)\\
\downarrow & & \downarrow & & \downarrow\\
H^0_{\et}(Y,\QQ_{\ell}) & \rightarrow & H^{2r}_{\et}(X, \QQ_{\ell}(r)) & \rightarrow &
H^{2r}_{\et}(Y, \QQ_{\ell}(r))
\end{array}
$$
in which the vertical arrows associate cohomology classes to cycles.  Note
that the composition of the two bottom maps is the map
$H^0_{\et}(Y,\QQ_{\ell}) \rightarrow H^{2r}_{\et}(Y,\QQ_{\ell}(r))$ induced by
$\eta_{\QQ_{\ell}}$.  Since $[Y] \in A^0(Y)$ maps to $1 \in H^0_{\et}(Y,\QQ_{\ell})$,
the commutativity of the above diagram implies that
$\eta_{\QQ_{\ell}}(1)$ is the cohomology class associated to $\pi^!\pi_*[Y]$,
as claimed.
\end{proof}

\begin{remark} \rm If $\pi$ is a closed immersion, the class $\pi^!\pi_*[Y]$
is simply the self-intersection of $Y$ in $X$, considered as a
cycle of codimension $r$ on $Y$.
\end{remark}

We now return to the situation considered at the beginning of this section.
The map ${\hat \pi}^{\prime}$ induces a map
$$H^j_{\et}(X_{U^{\prime}},\CL_{\xi^{\prime}}) \rightarrow 
H^j_{\et}(\hX_{U,\tau_0,i}^{\red}, ({\hat \pi}^{\prime})^*\CL_{\xi^{\prime}}).$$
Composing this with the map
$$\nu_{\hat \pi}: H^j_{\et}(\hX_{U,\tau_0,i}^{\red}, {\hat \pi}^*\CL_{\xi})
\rightarrow H^{j+2r}_{\et}(X_U,\CL_{\xi}(r)),$$
we obtain maps:
$$H^j_{\et}(X_{U^{\prime}},\CL_{\xi^{\prime}}) \rightarrow
H^{j+2r}_{et}(X_U,\CL_{\xi}(r)).$$
It is easy to verify that these maps are compatible with the action
of prime-to-$p$ Hecke operators on these cohomology spaces. 

Our main result is then:
\begin{theorem} \label{thm:main}
The maps:
$$H^j_{\et}(X_{U^{\prime}},\CL_{\xi^{\prime}}) \rightarrow
H^{j+2r}_{et}(X_U,\CL_{\xi}(r))$$
are injective.
\end{theorem}

The proof of this will occupy the remainder of this section, and the next.

Consider the class $c_{\hat \pi} \in 
H^{2r}_{\et}(\hX_{U,\tau_0,i}^{\red}, \QQ_{\ell}(r))$.  The map
$\hX_{U,\tau_0,i}^{\red} \rightarrow (X_{U^{\prime}})^{\tau_0,i}$
constructed in the previous section allows us to view this as a
cohomology class on $(X_{U^{\prime}})^{\tau_0,i}$.

The Leray spectral sequence for the natural map
$$\pi^{\prime}: (X_{U^{\prime}})^{\tau_0,i} \rightarrow X_{U^{\prime}}$$ 
is a spectral sequence
$$E_2^{j,k} = H^j_{\et}(X_{U^{\prime}}, R^k\pi^{\prime}_* \QQ_{\ell}) 
\rightarrow H^{j+k}_{\et}((X_{U^{\prime}})^{\tau_0,i}, \QQ_{\ell}).$$
It degenerates at $E_2$ by weight considerations.
In particular this yields a surjection 
$$H^{2r}_{\et}((X_{U^{\prime}})^{\tau_0,i}, \QQ_{\ell}) \rightarrow
H^0_{\et}(X_{U^{\prime}}, R^{2r}\pi^{\prime}_*\QQ_{\ell}).$$
Denote this surjection by $\alpha$.

Let $V$ be the complement of the (possibly empty) cycle 
$(X_{U^{\prime}})_{\tau_0,1}$, and let $V^{\tau_0,i}$ be the
preimage of $V$ in $X^{\tau_0,i}$.  Then $V^{\tau_0,i}$ is
a Grassmannian bundle over $V$, with fibers isomorphic to
$G(i,2i + r_{\tau_0} - r_{\sigma^{-1}\tau_0})$.  These fibers
have dimension $r$. 

Consider the map 
$$R^{2r}\pi^{\prime}_* \QQ_{\ell} \rightarrow 
j_*j^* R^{2r}\pi^{\prime}_* \QQ_{\ell},$$
where $j$ is the inclusion of $V$ in $X_{U^{\prime}}$.
Note that by the proper base change theorem, the stalk of 
$j^* R^{2r}\pi^{\prime}_* \QQ_{\ell}$ at a point $x$ of $V$ is 
isomorphic to $H^{2r}_{\et}(Z_x, \QQ_{\ell})$, 
where $Z_x = (\pi^{\prime})^{-1}(x)$, and is
therefore one-dimensional.  The map
$$H^{2r}_{\et}((X_{U^{\prime}})^{\tau_0,i}, \QQ_{\ell}) \rightarrow
H^{2r}_{\et}(Z_x, \QQ_{\ell})$$ 
given by applying $\alpha$
and then passing to the stalk at $x$ is the same as the map
induced by the inclusion of $Z$
in $(X_{U^{\prime}})^{\tau_0,i}$.

Let $W$ denote the universal $(\tau_0,i)$-constrained bundle on
$(X_{U^{\prime}})^{\tau_0,i}$.  It has rank $i$, and its restriction
$W_x$ to $Z_x$ for any point $x$ of $V$ can be
identified with the tautological
subbundle on $G(i,2i + r_{\tau_0} - r_{\sigma^{-1}\tau_0})$.

Denote by $c_i(W)$ the top Chern class of $W$, and consider the
class $C$ defined by
$$C = (-1)^rc_i(W)^{i + r_{\tau_0} - r_{\sigma^{-1}\tau_0}}.$$
For $x$ in $V$, the intersection $C \cap Z_x$
is $(-1)^rc_i(W_x)^{i + r_{\tau_0} - r_{\sigma^{-1}\tau_0}}$, which
is the class of a point on $Z_x$. 

Consider the class in 
$H^{2r}_{\et}((X_{U^{\prime}})^{\tau_0,i}, \QQ_{\ell}(r))$
arising from the cycle class $C$.  For each
$x$ in $V$, its restriction to
$H^{2r}_{\et}(Z_x, \QQ_{\ell}(r))$ is the fundamental
class of $H^{2r}_{\et}(Z_x, \QQ_{\ell}(r))$.  Thus the image 
under $\alpha$ of this class is an element of
$H^0_{\et}(X_{U^{\prime}}, R^{2r}\pi^{\prime}_*\QQ_{\ell}(r))$ that generates 
each stalk of $j^*R^{2r}\pi^{\prime}_*\QQ_{\ell}(r)$ as a $\QQ_{\ell}$-vector 
space.  
It follows that $j_*j^*R^{2r}\pi^{\prime}_*\QQ_{\ell}(r)$ is isomorphic to the 
constant sheaf $\QQ_{\ell}$, and that the map
$$R^{2r}\pi^{\prime}_*\QQ_{\ell}(r) \rightarrow 
j_*j^*R^{2r}\pi^{\prime}_*\QQ_{\ell}(r) \cong \QQ_{\ell}$$
is {\em split}.  In particular we obtain, for each $j$,
a surjection 
$$\beta_j: H^j_{\et}(X_{U^{\prime}}, R^{2r}\pi^{\prime}_*\QQ_{\ell}) \rightarrow
H^j_{\et}(X_{U^{\prime}}, \QQ_{\ell}(-r)).$$

Also note that by the projection formula, 
$$R^{2r}\pi^{\prime}_*(\pi^{\prime})^*\CL_{\xi^{\prime}} \cong
R^{2r}\pi^{\prime}_*\QQ_{\ell} \otimes_{\QQ_{\ell}} \CL_{\xi^{\prime}},$$
and therefore $\CL_{\xi^{\prime}}(-r)$ is a direct summand
of $R^{2r}\pi^{\prime}_*(\pi^{\prime})^*\CL_{\xi^{\prime}}$.
We therefore obtain for each $j$ a surjection
$$\beta_{j,\xi^{\prime}}:
H^p_{\et}(X_{U^{\prime}}, 
R^{2r}\pi^{\prime}_*(\pi^{\prime})^*\CL_{\xi^{\prime}}) 
\rightarrow H^j_{\et}(X_{U^{\prime}}, \CL_{\xi^{\prime}}(-r)).$$

The Leray spectral sequence induces an increasing filtration 
$\Fil^m_j(\CL_{\xi^{\prime}})$ on
$H^j_{\et}((X_{U^{\prime}})^{\tau_0,i},(\pi^{\prime})^*\CL_{\xi^{\prime}})$, 
such that 
$$\Fil^{m}_j(\CL_{\xi^{\prime}})/\Fil^{m-1}_j(\CL_{\xi^{\prime}}) =
H^m_{\et}(X_{U^{\prime}}, 
R^{j-m}\pi^{\prime}_*(\pi^{\prime})^*\CL_{\xi^{\prime}}).$$

Let $b$ be a class in $H^j_{\et}(X_{U^{\prime}}, \CL_{\xi^{\prime}})$, and
let $c$ be a class in $H^{2r}_{\et}(X_{U^{\prime}},\QQ_{\ell}(r))$.  Then
$(\pi^{\prime})^*b$ is a class in $\Fil^0_j$.  The product
$c \cup b$ then lies in $\Fil^{2r}_{j+2r}(\CL_{\xi^{\prime}}(r))$, and hence
maps onto $H^{j+2r}_{\et}(X_{U^{\prime}}, 
R^{2r}\pi^{\prime}_*(\pi^{\prime})^*\CL_{\xi^{\prime}}(r))$.
This in turn maps via $\beta_{j,\xi^{\prime}}$ onto 
$H^j_{\et}(X_{U^{\prime}}, \CL_{\xi^{\prime}})$.  

We thus obtain a map
from $H^j_{\et}(X_{U^{\prime}}, \CL_{\xi^{\prime}})$ to itself.  It can
be described in terms of $c$ in the following way: $\alpha(c)$ is
an element of $H^{2r}_{\et}(X_{U^{\prime}}, R^{2r}\pi^{\prime}_*\QQ_{\ell}(r))$;
this maps onto $H^0_{\et}(X_{U^{\prime}}, \QQ_{\ell})$ via $\beta_0$.  Thus
$\beta_0(\alpha(c))$ is a class of $H^0_{\et}(X_{U^{\prime}}, \QQ_{\ell})$;
the endomorphism of $H^j_{\et}(X_{U^{\prime}}, \CL_{\xi^{\prime}})$ described
above is simply multiplication by this class.

The upshot of all of this is:
\begin{proposition} \label{prop:cohomology}
If, for a particular choice of $X_U,\tau_0,$ and $i$, the corresponding
$\beta_0\alpha(c_{\hat \pi})$ is nonvanishing, then Theorem~\ref{thm:main}
holds for $X_U,\tau_0,$ and $i$ (and all $\CL_{\xi}$).
\end{proposition}
\begin{proof}
Consider the map
$$H^j_{\et}(X_{U^{\prime}}, \CL_{\xi^{\prime}}) \rightarrow
H^{j+2r}_{\et}((X_{U^{\prime}})^{\tau_0,i}, 
(\pi^{\prime})^*\CL_{\xi^{\prime}}(r))$$
that is obtained from the map
$$H^j_{\et}(X_{U^{\prime}}, \CL_{\xi^{\prime}}) \rightarrow
H^{j+2r}_{\et}(X_U, \CL_{\xi}(r))$$
of Theorem~\ref{thm:main} by composing with the natural map
$$H^{j+2r}_{\et}(X_U, \CL_{\xi}(r)) \rightarrow
H^{j+2r}_{\et}(\hX^{\red}_{U,\tau_0,i}, {\hat \pi}^*\CL_{\xi}(r)),$$
and identifying 
$H^{j+2r}_{\et}(\hX^{\red}_{U,\tau_0,i}, {\hat \pi}^*\CL_{\xi}(r))$
with $H^{j+2r}_{\et}((X_{U^{\prime}})^{\tau_0,i}, 
(\pi^{\prime})^*\CL_{\xi^{\prime}}(r)).$
To establish Theorem~\ref{thm:main}, it suffices to show this map is
injective.

By Lemma~\ref{lemma:reflect}, this map takes an element $b$
of $H^j_{\et}(X_{U^{\prime}},\CL_{\xi^{\prime}})$
to $c_{\hat \pi} \cup {\hat \pi}^*b$.  This lies in
$\Fil^{2r}_j(\CL_{\xi^{\prime}}(r))$, and maps via
$\beta_{j,\xi^{\prime}}$ to the element $\beta_0(\alpha(c_{\hat \pi})) b$
of $H^j_{\et}(X_{U^{\prime}},\CL_{\xi^{\prime}})$.  This element is
clearly nonzero if $b$ is.
\end{proof}

\section{The Thom-Porteus formula} \label{sec:porteus}

In this section we complete the proof of Theorem~\ref{thm:main} by
computing $\beta_0\alpha(c_{\hat \pi})$.  The key ingredient is the
Thom-Porteus formula, which will give us an expression for the
cycle class of $X_{U,\tau_0,i}$ (in the Chow ring of $X_U$) in
terms of a polynomial in Chern classes of bundles on $X_U$.

Before we state this formula we will need a bit of notation.  For $X$
a scheme, let $A^*(X) = \oplus_r A^r(X)$ denote the Chow ring of $X$.  
For an element $c$ in $A^*(X)$ we denote by $\Delta_q^{(p)}(c)$ the
determinant of the $p$ by $p$ matrix
$$\begin{pmatrix} 
c_q & c_{q+1} & \dots & c_{q + p - 1}\\
c_{q-1} & c_q & \dots & c_{q + p - 2}\\
\vdots & \vdots & & \vdots\\
c_{q - p + 1} & c_{q - p + 2} & \dots & c_q\\
\end{pmatrix}.
$$
Here $c_r$ is the $r$th graded part of $c$; the determinant $\Delta_q^{(p)}(c)$
therefore lies in $A^{pq}(X)$.

Then the Thom-Porteus formula states:
\begin{theorem}[\cite{Fulton}, 14.IV.4]
Let $X$ be a Cohen-Macaulay scheme, purely of dimension $N$, and
$\varsigma: E \rightarrow F$ a map of vector bundles on $X$.  Let
$D_k(\varsigma)$ denote the subscheme of $X$ defined by the condition
$\rank \varsigma \leq k$.  Suppose that $D_k(\varsigma)$ has the
``expected codimension''; that is, that the codimension of $D_k(\varsigma)$
is equal to $(e-k)(f-k)$ where $e$ and $f$ are the ranks of $E$ and $F$. 
Then, as elements of $A^{(e-k)(f-k)}(X)$, we have:
$$[D_k(\varsigma)] = \Delta_{f-k}^{(e-k)}(c(F)c(E)^{-1}).$$
\end{theorem}

We apply this to $X_{U,\tau_0,i}$.  For each $\tau$, let
$\CE_{\tau}$ denote the bundle $\Lie(A/X_U)^*_{p_{\tau}}$.  
\begin{proposition}
As elements of $A^r(X_U)$, we have:
$$[X_{U,\tau_0,i}] = 
\Delta_i^{(i+r_{\tau_0}-r_{\sigma^{-1}\tau_0})}
(c(\CE_{\tau_0})^{-1}
c(\Fabs^*\CE_{\sigma^{-1}\tau_0})),$$ 
where $\Fabs$ is the absolute Frobenius.
\end{proposition}
\begin{proof}
The cycle $X_{U,\tau_0,i}$ is the locus where the map
$$\Ver: \CE_{\tau_0} \rightarrow
\Lie(\CA^{(p)}/X_U)^*_{p_{\tau_0}}$$
has rank less than (or equal to) $r_{\sigma^{-1}\tau_0}-i$.
It has the expected codimension (equal to $r$), so the Thom-Porteus
formula yields:
$$[X_{U,\tau_0,i}] = 
\Delta_i^{(i+r_{\tau_0}-r_{\sigma^{-1}\tau_0})}
(c(\CE_{\tau_0})^{-1}
c(\Lie(\CA^{(p)}/X_U)^*_{p_{\tau_0}})).$$ 
The result then follows from the isomorphism:
$$\Lie(\CA^{(p)}/X_U)^*_{p_{\tau_0}} \cong 
\Fabs^*\Lie(\CA/X_U)^*_{p_{\sigma^{-1}\tau_0}} = \Fabs^*\CE_{\sigma^{-1}\tau}.$$
\end{proof}

This allows us to express ${\hat \pi}^!{\hat \pi}_*[\hX_{U,\tau_0,i}^{\red}]$
in terms of Chern classes of bundles on $\hX_{U,\tau_0,i}^{\red}$.
In particular we have:
$${\hat \pi}_*[\hX_{U,\tau_0,i}^{\red}] = [X_{U,\tau_0,i}].$$
Since for any bundle $F$ on $X_U$ we have ${\hat \pi}^!c(F) = c({\hat \pi}^*F)$
(see~\cite{Fulton}, Theorem 6.3 and the paragraph before example 8.1.1),
it follows that we have:
$${\hat \pi}^!{\hat \pi}_*[\hX_{U,\tau_0,i}^{\red}] =
\Delta_i^{(i+r_{\tau_0}-r_{\sigma^{-1}\tau_0})}
(c(\hCE_{\tau_0})^{-1}
c(\Fabs^*\hCE_{\sigma^{-1}\tau_0})),$$ 
where $\hCE_{\tau}$ is the restriction of $\CE_{\tau}$ to $\hX_{U,\tau_0,i}^{\red}$.

The next step is to express this in terms of classes of bundles
pulled back from $(X_{U^{\prime}})^{\tau_0,i}.$  Let 
$\iota$ denote the natural purely inseparable map 
$\hX_{U,\tau_0,i}^{\red} \rightarrow (X_{U^{\prime}})^{\tau_0,i},$
and recall that $W$ denotes the tautological (rank $i$) bundle
on $(X_{U^{\prime}})^{\tau_0,i}$.

Let $\CB$ be the universal abelian variety on $X_{U^{\prime}}$,
and let $\CE_{\tau}^{\prime}$ be the bundle $\Lie(\CB/X_{U^\prime})_{p_{\tau}}$.

\begin{proposition}
There are exact sequences:
$$0 \rightarrow \iota^*W \rightarrow 
({\hat \pi}^{\prime})^*\CE_{\tau_0}^{\prime} \rightarrow
\hCE_{\tau_0} \rightarrow 0$$
$$0 \rightarrow 
\Fabs^*({\hat \pi}^{\prime})^*\CE_{\sigma^{-1}\tau_0}^{\prime}
\rightarrow
\Fabs^*\hCE_{\sigma^{-1}\tau_0}
\rightarrow \iota^*W \rightarrow 0.$$
\end{proposition}
\begin{proof}
Let $\hCA$ be the restriction of $\CA$ to $\hX_{U,\tau_0,i}^{\red}$,
and $\hCB$ be the pullback of $\CB$ to $\hX_{U,\tau_0,i}^{\red}$.  The
universal isogeny $\hCA \rightarrow \hCB$ then induces a map
$$\Lie(\hCB/\hX_{U,\tau_0,i}^{\red})^*_{p_{\tau_0}} \rightarrow
\hCE_{\tau_0}.$$
The kernel of this map is $\iota^*W$; counting dimensions after pulling back
to any closed point shows it is surjective.
The first exact sequence above then follows
from the isomorphism 
$$\Lie(\hCB/\hX_{U,\tau_0,i}^{\red}) \cong 
(\hat \pi^{\prime})^*\CE_{\tau_0}^{\prime}.$$

For the second exact sequence, consider the commutative diagram:
$$\begin{array}{ccc}
\Fabs^*(H^1_{\DR}(\hCB/\hX_{U,\tau_0,i}^{\red})_{p_{\sigma^{-1}\tau_0}}) &
\rightarrow &
H^1_{\DR}(\hCB/\hX_{U,\tau_0,i}^{\red})_{p_{\tau_0}} \\
\downarrow & & \downarrow\\
\Fabs^*(H^1_{\DR}(\hCA/\hX_{U,\tau_0,i}^{\red})_{p_{\sigma^{-1}\tau_0}}) &
\rightarrow &
H^1_{\DR}(\hCA/\hX_{U,\tau_0,i}^{\red})_{p_{\tau_0}} \\
\end{array}$$
where the vertical maps are induced by the universal isogeny
$\hCA \rightarrow \hCB$, and the horizontal maps are induced by
relative Frobenius.

Note that the left-hand vertical map is an isomorphism.  Thus the
kernel of the bottom map (which is equal to 
$\Fabs^*\hCE_{\sigma^{-1}\tau_0}$) is 
isomorphic to the kernel of the composition of the upper horizontal map
and the right-hand vertical map.  The kernel of the former is
$\Fabs^*\Lie(\hCB/\hX_{U,\tau_0,i}^{\red})^*_{p_{\sigma^{-1}\tau_0}}$;
the kernel of the latter is $\iota^*W$.  We thus obtain an exact sequence:
$$0 \rightarrow
\Fabs^*\Lie(\hCB/\hX_{U,\tau_0,i}^{\red})^*_{p_{\sigma^{-1}\tau_0}}
\rightarrow
\Fabs^*\hCE_{\sigma^{-1}\tau_0} 
\rightarrow \iota^*W \rightarrow 0.$$
(Exactness on the right follows by pulling back to closed points and
counting dimensions.)  The second exact sequence again follows from
the isomorphism
$$\Lie(\hCB/\hX_{U,\tau_0,i}^{\red}) \cong 
({\hat \pi^{\prime}})^*\Lie(\CB/X_{U^{\prime}}).$$
\end{proof}

By the multiplicativity of the total Chern class, we have:
$${\hat \pi}^!{\hat \pi}_*[\hX_{U,\tau_0,i}^{\red}] =
\Delta_i^{(i+r_{\tau_0}-r_{\sigma^{-1}\tau_0})}
(c(\iota^*W)^2c(({\hat \pi^{\prime}})^*\CE^{\prime}_{\tau_0})^{-1}
c(\Fabs^*({\hat \pi^{\prime}})^*\CE^{\prime}_{\sigma^{-1}\tau_0})).$$ 
This means that $c_{\hat \pi}$, when considered as a cohomology class
on $(X_{U^{\prime}})^{\tau_0,i}$, is the cohomology class associated to
the element:
$$\Delta_i^{(i+r_{\tau_0}-r_{\sigma^{-1}\tau_0})}
(c(W)^2c((\pi^{\prime})^*\CE^{\prime}_{\tau_0})^{-1}
c((\pi^{\prime})^*\Fabs^*\CE^{\prime}_{\sigma^{-1}\tau_0}))
\in A^r((X_{U^{\prime}})^{\tau_0,i}).$$ 

Consider $\beta_0(\alpha(c_{\hat \pi}))$.  It is a section of
the constant sheaf $\QQ_{\ell}$ on $(X_{U^{\prime}})$.  We have shown that
for any $x$ in $V$, $\beta_0(\alpha(c_{\hat \pi}))_x$ is obtained
by pulling back $c_{\hat \pi}$ to an element  of 
$H^{2r}_{\et}((\pi^{\prime})^{-1}(x), \QQ_{\ell}(r))$, and applying the 
canonical 
isomorphism of this space with $\QQ_{\ell}$.  In other words, 
$\beta_0(\alpha(c_{\hat \pi}))_x$
is the element of $H^{2r}_{\et}((\pi^{\prime})^{-1}(x), \QQ_{\ell}(r))$
obtained by intersecting 
$$\Delta_i^{(i+r_{\tau_0}-r_{\sigma^{-1}\tau_0})}
(c(W)^2c((\pi^{\prime})^*\CE^{\prime}_{\tau_0})^{-1}
c((\pi^{\prime})^*\Fabs^*\CE^{\prime}_{\sigma^{-1}\tau_0}))$$ 
with $Z = (\pi^{\prime})^{-1}(x)$ and then taking the associated cohomology 
class.

Any bundle pulled back from $X_{U^{\prime}}$ via $\pi^{\prime}$ becomes
the trivial bundle when restricted to $Z$, and
thus has trivial total Chern class in $A^*(Z).$
On the other hand, $Z$ is a Grassmannian of
$i$ planes in a $2i + r_{\tau_0} - r_{\sigma^{-1}\tau_0}$ dimensional
space, and $W$ restricts to the tautological subbundle $W_Z$
on this Grassmannian.

Thus, $\beta_0(\alpha(c_{\hat \pi}))$ is the cohomology class in
$H^{2r}_{\et}(Z, \QQ_{\ell}(r))$ associated to
the element 
$$\Delta_i^{(i+r_{\tau_0}-r_{\sigma^{-1}\tau_0})}(c(W_Z)^2) \in
A^r(Z).$$

\begin{proposition}
Let $[P]$ be the class of a point in $A^r(Z)$.
Then we have:
$$\Delta_i^{(i+r_{\tau_0}-r_{\sigma^{-1}\tau_0})}(c(W_Z)^2) = 
(-1)^r\binom{2i + r_{\tau_0} - r_{\sigma^{-1}\tau_0}}{i}[P].$$
\end{proposition}
\begin{proof}
We have a tautological exact sequence:
$$0 \rightarrow W_Z \rightarrow 
\OO_Z^{2i + r_{\tau_0} - r_{\sigma^{-1}\tau_0}}
\rightarrow Q \rightarrow 0$$
of vector bundles on $Z$.
Dualizing yields a sequence:
$$0 \rightarrow Q^* \rightarrow 
\OO_Z^{2i + r_{\tau_0} - r_{\sigma^{-1}\tau_0}}
\rightarrow W_Z^* \rightarrow 0.$$

Let $M$ be any endomorphism of 
$\overline{\FF}_p^{2i + r_{\tau_0} - r_{\sigma^{-1}\tau_0}}$
with distinct eigenvalues.  We obtain a map 
$\varsigma_M: Q^* \rightarrow W_Z^*$ by including $Q^*$ in
$\OO_Z^{2i + r_{\tau_0} - r_{\sigma^{-1}\tau_0}}$,
applying the endomorphism $M$, and then projecting to $W_Z^*$.

The subscheme $D_0(\varsigma_M)$ of points of $Z$ on 
which $\varsigma_M$ is the zero map is easily seen to be reduced, and equal
to the union of those points of $Z$ that correspond
to $i + r_{\tau_0} - r_{\sigma^{-1}\tau_0}$-dimensional subspaces
of $\overline{\FF}_p^{2i + r_{\tau_0} - r_{\sigma^{-1}\tau_0}}$
that are stable under $M$.  Any such space is the direct sum of precisely
$i + r_{\tau_0} - r_{\sigma^{-1}\tau_0}$ of the
$2i + r_{\tau_0} - r_{\sigma^{-1}\tau_0}$ distinct eigenspaces of $M$.
Thus we have:
$$[D_0(\varsigma_M)] = 
\binom{2i + r_{\tau_0} - r_{\sigma^{-i}\tau_0}}{i}[P].$$

On the other hand, the Thom-Porteus formula tells us that we have:
$$[D_0(\varsigma_M)] = 
\Delta_i^{(i + r_{\tau_0} - r_{\sigma^{-1}\tau_0})}(c(W_Z^*)c(Q^*)^{-1}).$$
The result follows by putting these two together, and using the basic 
identities:
$$c(W_Z^*)c(Q^*) = 1,$$
$$c_j(W_Z) = (-1)^jc_j(W_Z^*).$$
\end{proof}

It follows that $\beta_0(\alpha(c_{\hat \pi}))$ is a non-vanishing
section of the constant sheaf $\QQ_{\ell}$ on $X_{U^{\prime}}.$  
Theorem~\ref{thm:main} thus follows from Proposition~\ref{prop:cohomology}.

\section{Jacquet-Langlands correspondences} \label{sec:JL}

We now use the above characteristic $p$ results to study the cohomology
of Shimura varieties in characteristic zero.  As the Shimura varieties
we consider are not necessarily proper, we will first need some results beyond
the standard theory of vanishing cycles to accomplish this. 

Let $S = \spec W(\overline{\FF}_p)$; let $s$ and $\eta$ denote the closed
point and a geometric generic point of $S$, respectively.  Let $X$
be a smooth $S$-scheme, and let $\overline{X}$ be a compactification of
$X$ with the following properties:

\begin{itemize}
\item $\overline{X}$ is smooth and proper over $S$,
\item the complement $\overline{X}\setminus X$ is a divisor $D$ with normal
crossings, and
\item if $D_1, \dots, D_i$ are irreducible components of $D$, then the
intersection $D_1 \cap \dots \cap D_i$ is either empty or smooth over $S$.
\end{itemize}

Under these hypotheses, we have:
\begin{lemma}
Let $\overline{\CF}$ be a lisse sheaf on $\overline{X}$, and let
$\CF$ denote its restriction to $X$.  Then the specialization maps
$$H^j_{\et}(X_{\eta},\CF_{\eta}) \rightarrow H^j_{\et}(X_s,\CF_s)$$
are isomorphisms.
\end{lemma}
\begin{proof}
We work by induction on the number of irreducible components of $D$.
In the base case $D$ is empty and the above result is immediate from the
theory of vanishing cycles.

Suppose the result is true for $D$ having $k$ components.  Let 
$X^k = \overline{X} \setminus (D_1 \cup \dots \cup D_k)$,
and let $D_{k+1}^k = D_{k+1} \setminus (D_1 \cup \dots \cup D_k)$.
Then the
specialization maps:
$$H^j_{\et}(X^k_{\eta},\CF_{\eta}) \rightarrow H^j_{\et}(X^k_s, \CF_s)$$
$$H^j_{\et}((D_{k+1}^k)_{\eta},\CF_{\eta}) \rightarrow 
H^j_{et}((D_{k+1}^k)_s, \CF_s)$$
are isomorphisms.

These fit into a Gysin sequence:
$$
\begin{array}{ccccccc}
\rightarrow & H^j_{\et}((D_{k+1}^k)_{\eta},\CF_{\eta}) 
      & \rightarrow & H^{j+2}_{\et}(X^k_{\eta},\CF_{\eta})
      & \rightarrow & H^{j+2}_{\et}(X^{k+1}_{\eta},\CF_{\eta})
      & \rightarrow\\
& \downarrow & & \downarrow & & \downarrow & \\
\rightarrow & H^j_{\et}((D_{k+1}^k)_s,\CF_s) 
      & \rightarrow & H^{j+2}_{\et}(X^k_s,\CF_s)
      & \rightarrow & H^{j+2}_{\et}(X^{k+1}_s,\CF_s)
      & \rightarrow\\
\end{array}
$$
and hence the result holds for $X^{k+1}$ as well.
\end{proof}

Call such a compactification of $X$ a {\em good} compactification.
In \cite{FC}, Faltings-Chai show that the toroidal compactifications
of the moduli spaces of principally polarized abelian varieties are
good compactifications, and assert (with no details) that their methods
carry over to arbitrary PEL Shimura varieties.  As yet unpublished
work of Kai-Wen Lan (\cite{La}, Thm 6.4.1.1) shows that the unitary 
Shimura varieties $X_U$
admit good (toroidal) compactifications $\overline{X}_U$.  Moreover,
the sheaves
$\CL_{\xi}$ extend to lisse sheaves on $\overline{X}_U$.  We thus have
a natural isomorphism:
$$H^j_{\et}((X_U)_{\eta},(\CL_{\xi})_{\eta}) \cong
H^j_{\et}((X_U)_s,(\CL_{\xi})_s).$$

We can use this, together with Theorem~\ref{thm:main}, to ``transfer''
automorphic representations from one algebraic group to another.
Fix two places $\tau_0$, $\tau_0^{\prime}$ of $F^+$, with
$r_{\tau_0^{\prime}} \leq r_{\tau_0}$, and fix an $i$
such that $1 \leq i \leq \min(r_{\tau_0^{\prime}},n-r_{\tau_0})$.
Then there exists a $\CV^{\prime}$ such that $\CV^{\prime}(\AA^{\infty}_{\QQ})$
is isomorphic to $\CV(\AA^{\infty}_{\QQ})$, but whose invariants at infinity
satisfy:
\begin{itemize}
\item $r_{\tau_0}(\CV^{\prime}) = r_{\tau_0} + i$,
\item $r_{\tau_0^{\prime}}(\CV^{\prime}) = r_{\tau_0^{\prime}} - i$,
\item $r_{\tau}(\CV^{\prime}) = r_{\tau}(\CV)$ for $\tau$ outside
$\{\tau_0,\tau_0^{\prime}\}$.
\end{itemize}
Fix such a $\CV^{\prime}$, and let $G^{\prime}$ be the corresponding 
unitary group.  Also fix an identification of $\CV(\AA^{\infty}_{\QQ})$
with $\CV^{\prime}(\AA^{\infty}_{\QQ})$; this yields an
identification of $G(\AA^{\infty}_{\QQ})$ with $G^{\prime}(A^{\infty}_{\QQ})$.

\begin{theorem} \label{thm:JL}
Let $\pi^{\prime}$ be an automorphic representation
of $G^{\prime}$, and suppose that there exists a representation
$\xi^{\prime}$ of $G^{\prime}(\AA^{\infty}_{\QQ})$ over $\overline{\QQ}_{\ell}$
such that $\pi^{\prime}_{\infty}$ is cohomological for $\xi^{\prime}$.
Suppose also that there exist good compactifications for 
unitary Shimura varieties attached to $G$ and $G^{\prime}$.
Then there exists an automorphic representation $\pi$ of $G$ such
that $\pi_v = \pi^{\prime}_v$ for all finite places $v$ of $\QQ$,
and such that $\pi_{\infty}$ is cohomological for the
representation $\xi$ of $G(\AA^{\infty}_{\QQ})$ that corresponds to
$\xi^{\prime}$.
\end{theorem}

\begin{proof}
Let $U^{\prime}$ be a compact open subgroup of 
$G^{\prime}(\AA^{\infty}_{\QQ})$, such that $\pi^{\prime}$ has a
nonzero $U^{\prime}$-fixed vector.  Let $U$ be the corresponding
subgroup of $G(\AA^{\infty}_{\QQ})$.

Fix, for each $p$, an embedding of $W(\overline{\FF}_p)$ as a subring of
$\CC$.  This determines a Frobenius element $\Frob_p$ of 
$\gal(\overline{\QQ}/\QQ)$, up to inertia.  By {\v C}ebotarev,
we can find a $p$ such that $p$ is unramified in $F$ and split in $E$,
such that $U_p$ is a maximal compact subgroup of $G(\QQ_p)$, and such that
$\Frob_p\tau_0^{\prime} = \tau_0$.  Also choose an auxiliary prime $l$
different from $p$.  

Associated to these choices we have Shimura varieties $X_U$ and
$X_{U^{\prime}}$ over $W(\overline{\FF}_p)$.  Let $N$ be
an integer divisible by all the primes of bad reduction of $X_U$,
and let $\TT_U$ be the Hecke algebra (over $\QQ_{\ell}$)
of prime-to-$Np$ Hecke operators for $G$.

Let $s: \spec \overline{\FF}_p
\rightarrow \spec W(\overline{\FF}_p)$ be the closed point
of $\spec W(\overline{\FF}_p)$.  Theorem~\ref{thm:main}
yields, for all $j$, an injection
$$H^j_{\et}((X_{U^{\prime}})_s, (\CL_{\xi^{\prime}})_s) \rightarrow
H^{j+2r}_{\et}((X_U)_s, (\CL_{\xi})_s)$$
that is compatible with the action of the Hecke algebra $\TT_U$. 
This yields $\TT_U$-equivariant
injections
$$H^j_{\et}((X_{U^{\prime}})_{\eta}, (\CL_{\xi^{\prime}})_{\eta}) \rightarrow
H^{j+2r}_{\et}((X_U)_{\eta}, (\CL_{\xi})_{\eta}),$$ 
by the above lemma, where $\eta$ is a geometric generic point of
$\spec W(\overline{\FF}_p)$.

The representation $\pi^{\prime}$ determines a maximal ideal $\mm$ of
$\TT_U$, and our hypotheses guarantee that for some $j$,
$H^j_{\et}((X_{U^{\prime}})_{\eta}, (\CL_{\xi^{\prime}})_{\eta})_{\mm}$
will be nonzero.  Then
$H^{j+2r}_{\et}((X_U)_{\eta}, (\CL_{\xi})_{\eta})_{\mm}$ is nonzero as well.
There is therefore an automorphic representation $\pi$ of $G$,
such that: 
\begin{itemize}
\item $\pi_{\infty}$ is cohomological for $\xi$, 
\item $\pi$ has a $U$-fixed vector, and 
\item for any Hecke operator in $\TT_U$, the Hecke eigenvalue for $\pi$ is
the same as that for $\pi^{\prime}$.
\end{itemize}

It follows that $\pi_v$ is isomorphic to $\pi^{\prime}_v$ for any finite
place $v$ not dividing $Np$.  By Cebotarev, the Galois representations
associated to $\pi$ and $\pi^{\prime}$ coincide.  It then follows that
$\pi_v = \pi^{\prime}_v$ for all finite places $v$.
\end{proof}

\textsc{Acknowledgements}

The author is deeply indebted to Richard Taylor for his continued interest
and encouragement.  This research was partially supported by the National 
Science Foundation.

\end{document}